\documentclass[11pt,a4paper,twoside]{amsart}
\usepackage{amsmath}
\usepackage{amssymb}
\usepackage{latexsym}

\newcommand{\co}{\colon\thinspace}    
\newcommand{\fnote}[1]{\footnote{\small sharp1}}
\newcommand{\inv}{^{-1}}              

\newcommand{\N}{{\mathbb N}}
\newcommand{\Z}{{\mathbb Z}}
\newcommand{\R}{{\mathbb R}}

\newcommand{\T}{{\mathbb T}}
\newcommand{\Sbar}{{\mathbb S}}

\newcommand{\spt}{\mbox{supp}}
\newcommand{\Azero}{\mathcal{A}_0}

\newcommand{\tildeAzero}{\tilde{\mathcal{A}_0}}
\newcommand{\inter}{\mbox{Int}}

\newtheorem{theorem}{Theorem}[section]
\newtheorem{proposition}[theorem]{Proposition}
\newtheorem{corollary}[theorem]{Corollary}

\newtheorem{lemma}[theorem]{Lemma}
\newtheorem{remark}[theorem]{Remark}

\title{Aubry sets vs Mather sets \\
in two degrees of freedom }
\author{Daniel Massart}
\date{\today}
\begin{document}

\begin{abstract}
 Let $L$ be an autonomous Tonelli Lagrangian on a closed manifold of dimension two. Let $\mathcal{C}$ be the set of   cohomology classes whose Mather set consists of periodic orbits, none of which is a fixed point. Then for almost all $c$ in $\mathcal{C}$, the Aubry set of $c$ equals the Mather set of $c$.

\end{abstract}
\maketitle
\section{Introduction}
\subsection{Motivation}
We study Tonelli Lagrangian systems on closed manifolds, along the lines of \cite{Mather91}.  
The Aubry set is a specific invariant set of the Euler-Lagrange flow, originally defined in \cite{Mather93}, although its current name comes from \cite{Fathi}. Roughly speaking, it is the obstruction to push a Lagrangian submanifold inside a convex hypersurface of the cotangent bundle of a closed manifold  without changing its cohomology class (see \cite{PPS}). 
Various nice results hold when the Aubry set is a finite union of hyperbolic, periodic orbits : 
\begin{itemize}
  \item asymptotic estimates for near-optimal periodic geodesics (\cite{A03}), if the Lagrangian is a metric of negative curvature on a surface
  \item existence of "physical"  solutions of the Hamilton-Jacobi equation (\cite{AIPS05})
	\item existence of $C^{\infty}$ subsolutions of the Hamilton-Jacobi equation (\cite{Bernard}).
\end{itemize}

By \cite{CI} when there is a minimizing periodic orbit, a small perturbation makes it hyperbolic while still minimizing. The trouble is to find minimizing periodic orbits.

While this seems out of reach for the time being, there is a particular case where this difficulty is easily overcome: that is when the dimension of the configuration space is two, for then Proposition 2.1 of \cite{CMP}  says that any minimizing measure with a rational homology class is supported on periodic orbits. 

Even then, yet another problem arises: the Aubry set always contains the union of the supports of all minimizing measures (Mather set), but the inclusion may be proper. The purpose of this paper is to clarify the relationship between the Aubry set and the Mather set, when the latter consists of periodic orbits. In loose terms our main result says that in that case (and in two degrees of freedom) they almost always coincide. See the next paragraph. This is a generalization of a  result of Mather for twist maps, see \cite{Mather08}, section 3.

\textbf{Acknowledgements: } I thank Albert Fathi for his careful reading of the manuscript. This work was partially supported by the ANR project "Hamilton-Jacobi et th\'eorie KAM faible".
\subsection{Definitions and precise statements}
Recall that an  autonomous Tonelli Lagrangian on a closed manifold $M$ is  a $C^2$ function $L$ from $ TM$ to $ \R$ which is fiberwise  superlinear and such that $\partial^2 L / \partial v^2$ is positive definite everywhere. Let
\begin{itemize}
	\item $L$ be an autonomous Tonelli  Lagrangian on a closed manifold $M$ 
	\item  $\phi_t$ be the Euler-Lagrange flow of $L$
	\item $p$ be the canonical projection $TM \longrightarrow M$.
\end{itemize}

The first object one encounters when using variational methods is Ma\~n\'e's action potential  : for each nonnegative $t$, and $x,y$ in $M$, define
	\[ h_t (x,y) := \inf \int^{t}_{0}L\left(\gamma (s),\dot{\gamma}(s)\right) ds
\]
over all absolutely continuous curves $\gamma : \left[0,t\right]\longrightarrow M$ such that $\gamma (0)=x$, $\gamma (t)=y$. The infimum is in fact a minimum due to the fiberwise strict convexity and superlinearity  of $L$, and the curves achieving the minimum are projections to $M$ of pieces of orbits of $\phi_t$. Such curves are called extremal. 

Looking for orbits that realize the action potential between any two of their points, one is led to consider the Peierls barrier (\cite{Mather93})
	\[h (x,y) := \liminf_{t \rightarrow \infty}h_t (x,y).
\]
The projected Aubry set is then defined  as
	\[
	\mathcal{A}(L) := \left\{ x\in M \co h(x,x)=0\right\}.
\]
 Mather's Graph Theorem (\cite{Mather91}, see also \cite{Fathi}, Theorem 5.2.8) then says that for any $x \in \mathcal{A}(L)$, there exists a unique $v \in T_x M$ such that $p \circ \phi_t (x,v)$,  $t \in \R$, is an extremal curve. The set 
	\[
	\tilde{\mathcal{A}}(L) := \left\{ (x,v)\in M \co p \circ \phi_t (x,v) \in \mathcal{A}(L) \  \forall t \in \R \right\}
\]
is called the Aubry set of $L$, it is compact and $\phi_t$-invariant.

As noticed by Mather, it is often convenient to deal with invariant measures rather than individual orbits. Define $\mathcal{M}_{inv}$
to be the set of $\Phi_t$-invariant, compactly supported, Borel probability measures on $TM$.
Mather showed that the function (called action of the Lagrangian on measures)
	\[
	\begin{array}{rcl}
\mathcal{M}_{inv} & \longrightarrow & \R \\
\mu & \longmapsto & \int_{TM}	L d\mu
\end{array}
\]
is well defined and has a minimum.  A measure achieving the minimum is called $L$-minimizing. The value of the minimum, times minus one,  is called the critical value of $L$, and denoted $\alpha (L)$. It is also called effective Hamiltonian, see for instance \cite{Evans}. The Mather set $\tilde{\mathcal{M}}(L)$ of $L$ is then defined as the closure of the union of the supports of all minimizing measures. It is compact, $\phi_t$-invariant, and contained in $\tilde{\mathcal{A}}(L)$. 

The minimization procedure may be refined as follows. Mather observed that if $\omega$ is a closed one-form on $M$ and $\mu \in \mathcal{M}_{inv}$ then the integral $\int_{TM}	\omega d\mu$ is well defined, and only depends on the cohomology class of $\omega$. By duality this endows  $\mu$ with a homology class : $\left[\mu\right]$ is the unique $h \in H_1 (M,\R)$ such that  
	\[
\langle h,\left[\omega \right]\rangle = \int_{TM}	\omega d\mu 
\]
for any closed one-form $\omega$ on $M$. Besides, for any $h \in H_1 (M,\R)$, the set 
	\[ \mathcal{M}_{h,inv}:= \left\{\mu \in \mathcal{M}_{inv} \co \left[\mu\right]=h\right\}
\]
is not empty. Again the action of the Lagrangian on this smaller set of measures has a minimum, which is a function of $h$, called the $\beta$-function of the system, or effective Lagrangian. A measure achieving the minimum is called $(L,h)$-minimizing, or $h$-minimizing for short.  

When the dimension of $M$ is two, we get a bit of help from the topology. Let 
$\Gamma$ be the quotient of $H_1 (M,\Z)$ by its torsion (we do not assume $M$ to be orientable), $\Gamma$ embeds as a lattice into $H_1 (M,\R)$. A homology class $h$ is said to be 1-irrational if there exist $h_0 \in \Gamma$ and $r \in \R$ such that $h=rh_0$.
Proposition 2.1 of \cite{CMP} (see also Proposition 5.6 of \cite{nonor}) reads  :
\begin{proposition}\label{rational}
Let $M$ be a closed surface, possibly non-orientable, and let $L$ be a Tonelli Lagrangian on $M$.
If $h$ is a 1-irrational homology class and $\mu$ is an $h$-minimizing measure, then the support of $\mu$ consists of periodic orbits, or fixed points.
\end{proposition}

There is a dual construction : if $\omega$ is a closed one-form on  $M$, then $L-\omega$ is a Tonelli Lagrangian, and furthermore $L-\omega$ has the same Euler-Lagrange flow as $L$. The Aubry set, Mather set, and critical value of $L-\omega$ are denoted $\tilde{\mathcal{A}}_L (c), \tilde{\mathcal{M}}_L (c), \alpha_L (c)$ respectively, or just $\tilde{\mathcal{A}}(c), \tilde{\mathcal{M}}(c), \alpha (c)$ when no ambiguity is possible.  An  $(L-\omega)$-minimizing measure is also called $(L,\omega)$-minimizing, $(L,c)$-minimizing, or just $c$-minimizing for short if $c$ is the cohomology of $\omega$. In formal terms we have defined
	\[
	\begin{array}{rcl}
\beta_L \co H_1 (M,\R) & \longrightarrow & \R \\
h & \longmapsto & 
\min \left\{\int_{TM}	Ld\mu \co \left[\mu\right]=h\right\}\\
\alpha_L \co H^1 (M,\R) & \longrightarrow & \R \\
c & \longmapsto & 
-\min \left\{\int_{TM}	(L-\omega)d\mu \co \left[\omega\right]=c\right\}.
\end{array}
\]

Mather proved that $\alpha_L$ and $\beta_L$ are convex, superlinear, and Fenchel dual of one another. In particular $\min \alpha = -\beta(0)$, and we have the Fenchel inequality :
	\[\alpha_L (c) + \beta_L (h) \geq \left\langle c,h \right\rangle \  \forall c \in H^1 (M,\R), \; h \in H_1 (M,\R).
\]
Given $c \in H^1 (M,\R)$ (resp. $h \in H_1 (M,\R)$), the set of $h \in H_1 (M,\R)$ (resp. $c \in H^1 (M,\R)$) achieving equality in the Fenchel equality is called the Legendre transform of $c$ (resp. $h$), and denoted $\mathcal{L}(c)$ (resp. $\mathcal{L}(h)$). 

The functions $\alpha_L$  and $\beta_L$ are sometimes called homogeneized Hamiltonian and Lagrangian, respectively.

The main geometric features of a convex function are its smoothness and strict convexity, or lack thereof.  In general, the maps $\alpha_L$ and $\beta_L$ are neither strictly convex, nor smooth (\cite{gafa}). The regions where either map  is not strictly convex are called flats (see Appendix A for precise definitions). A flat  is a convex subset of a linear space, hence it makes sense to speak of its relative interior, or interior for short. The sets $\mathcal{L}(c)$ for $c \in H^1 (M,\R)$ (resp. $\mathcal{L}(h)$ for $h \in H_1 (M,\R)$), if they contain more than one point, are non-trivial instances of flats ; conversely, by the Hahn-Banach Theorem, any flat is contained in the Legendre transform of some point.

Note that if two cohomology classes lie in the relative interior of a flat $F$ of $\alpha_L$, by \cite{Mather91} their Mather sets coincide. We denote by $\tilde{\mathcal{M}}(F)$ the common Mather set to all the cohomologies in the relative interior of $F$. For any $c$ in $F$, the Mather set of $c$  contains the Mather set of $F$. We say a flat is rational if its Mather set consists of periodic orbits or fixed points. It is easy to see that any rational flat of $\alpha_L$ is contained in $\mathcal{L}(h)$ for some 1-irrational $h$. A partial converse is true when the dimension of $M$ is two (see Lemma \ref{nonsingular_3}).

 As to Aubry sets, Proposition 6 of \cite{ijm} reads :
 \begin{proposition} \label{pomface}
If a cohomology class $c_1$ belongs to 
a flat $F_{c}$ of $\alpha_L $ containing 
$c$ in its interior,
then $\mathcal{A}(c) \subset {\mathcal A}(c_1)$. 
In particular, if $c_1$ lies in 
the interior of $F_{c}$, then $\mathcal{A}(c) = {\mathcal A}(c_1)$. 
Conversely,
if two cohomology classes $c$ and $c_1$ are such that 
$\tilde{\mathcal{A}}(c) \cap \tilde{\mathcal {A}}(c_1) \neq \emptyset$, then $\alpha_L $ has a flat containing $c$ and $c_1$.
\end{proposition}
So for any flat $F$ of $\alpha$ and any $c_1$, $c_2$ in the interior of $F$, the Aubry sets $\tilde{\mathcal{A}}(c_1)$ and $\tilde{\mathcal{A}}(c_2)$ coincide.  We denote by $\tilde{\mathcal{A}}(F)$ the common Aubry set to all the cohomologies in the interior of $F$.

A flat of $\alpha_L $ is called singular if its Mather set  contains a fixed point of the Euler-Lagrange flow. A homology class $h$ is called singular if its Legendre transform is a singular flat. So the set of singular classes is either empty, or it contains zero and is compact.  When there are fixed points, we lose some of the perks of the low dimension, which explains why we have to exclude singular classes from our main result. The purpose of this paper is to prove that the Aubry set of a nonsingular rational flat equals its Mather set.
\begin{theorem}\label{principal}
Assume
\begin{itemize}
	\item $M$ is a closed surface
	\item $L$ is an autonomous  Tonelli Lagrangian on $M$
	\item $h$ is a 1-irrational, nonsingular homology class.
\end{itemize}
Then $\tilde{\mathcal{A}}(\mathcal{L}(h))=\tilde{\mathcal{M}}(\mathcal{L}(h))$, and $\tilde{\mathcal{A}}(\mathcal{L}(h))$ is a union of periodic orbits.
\end{theorem}
So in the interior of a nonsingular rational flat, the Aubry set is as small as possible since it must contain the Mather set.  Note that the boundary of a convex set $C$ is negligible in $C$, in any reasonable sense of negligible, which accounts for the phrase 'almost always coincide' we used in the 'Motivation' subsection.

Let us briefly review what is known in one degree of freedom. Take $L$ to be a time-periodic Lagrangian on the circle, and take a rational element  $h$ of the homology of the circle. Then $\mathcal{L}(h)$ is an interval $\left[c^-, c^+\right]$. It is stated in \cite{Mather93}, p. 1376, and proved in \cite {Mather08}, section 3, that for $c$ in the relative interior of $\left[c^-, c^+\right]$, then $\tilde{\mathcal{A}}(c)=\tilde{\mathcal{M}}(c)$, and 
$\tilde{\mathcal{A}}(c)$ consists of periodic orbits.

Here is an outline of this paper. In section 2 we prove a local result (Lemma \ref{local}) : if for some cohomology class $c$, the Aubry set $\mathcal{A}(c)$ contains a periodic orbit $\gamma$, then there exists a face $F$ of $\alpha$, containing $c$, although not necessarily in its interior,  such that the Aubry set $\mathcal{A}(F)$ (which is a subset of $\mathcal{A}(c)$ by Proposition \ref{pomface}), in a neighborhood of $\gamma$, contains only periodic orbits homotopic to $\gamma$.

From this we deduce (Corollary \ref{corollaire_local}) that if for some cohomology class $c$, the Mather set $\mathcal{M}(c)$ consists of periodic orbits, then there exists a face $F$ of $\alpha$, containing $c$, such that $\mathcal{A}(F)=\mathcal{M}(c)$. 

Our theorem would follow if we could show that for any 1-irrational, non-singular homology class $h$,  $\mathcal{M}(\mathcal{L}(h))$ consists of periodic orbits. This, by Lemma \ref{face_rationnelle}, would follow from $h$ lying in the relative interior of 
$\mathcal{F}(\mathcal{L}(h))$ (see the definition in Appendix A). However, this may not be true; but we prove in Section 3 that for some $h'$ in $R_h$, the largest radial flat of $\beta$ containing $h$, $h'$ is contained  in the relative interior of 
$\mathcal{F}(\mathcal{L}(h))$, and that is enough to prove Theorem \ref{principal}. 

In Section 4 we give some examples to show that the non-singularity hypothesis in our result is genuinely necessary. In Appendix A we have gathered the notions of convex analysis that we use. In Appendix B we prove the lemmas needed to include the case of non-orientable surfaces. In Appendix C we prove some technical results about the faces of $\beta$ (some of which were proved in \cite{gafa}, in the case when the Lagrangian is a Finsler metric) that we need for Lemma \ref{nonsingular_3}.

\section{Local structure of the Aubry set at periodic orbits}
Our next lemma is   a slight modification of  Proposition 5.4 of \cite{nonor}. In \cite{nonor} only geodesic flows are considered but the proof extends without modification to the case of Lagrangian flows.

We say a closed curve is minimizing (resp. $c$-minimizing) if the probability measure equidistributed on it is minimizing (resp. $c$-minimizing). This is equivalent to saying that the closed curve is an extremal contained in the projected Mather set.

\begin{lemma}\label{key_o}
Let \begin{itemize}
  \item $M$ be an oriented closed surface
  \item $L$ be an autonomous Tonelli Lagrangian on $M$
  \item $\gamma_0$ be a closed,  minimizing extremal of $L$, such that $\gamma_0$ is not a fixed point  
  \item $h_0$ be the homology class of the minimizing measure supported on 
  $(\gamma_0,\dot\gamma_0)$
  \item $c$ be a cohomology class in $\mathcal{L}(h_0)$.
    \end{itemize}
 There exists a neighborhood $V_0$ of $(\gamma_0,\dot\gamma_0)$ in $TM$ such that for any simple  extremal $\gamma$ such that $(\gamma,\dot\gamma)$ is contained in $\tilde{\mathcal{A}}(c)$, if  $(\gamma,\dot\gamma)$ enters (resp. leaves) $V_0$  then $\gamma$ is either 
 \begin{itemize}
  \item a closed extremal asymptotic to $\gamma_0$
  \item or positively (resp. negatively) asymptotic to a closed extremal asymptotic to $\gamma_0$.
 \end{itemize}
  \end{lemma}
\subsection{}
Let $\gamma_0$ be a $C^1$ simple closed curve (not a fixed point) in an oriented surface $M$. We say a $C^1$ curve $\alpha \co \R \longrightarrow M \setminus \gamma_0$ is positively asymptotic to $\gamma_0$ on the right if the $\omega$-limit set of $\alpha$ is $\gamma_0$ and there exists some $t_0 \in \R$ such that, for any $t \geq t_0$, $\alpha (t)$ lies in the right-hand  side (with respect to the chosen orientation of $M$) of a tubular neighborhood of $\gamma_0$. Similar definitions can be made replacing positively with negatively, and right with left. The lemma below will be used in the proof of Lemma \ref{local}.
\begin{lemma}\label{asymptotic}
Let $\gamma_0$ be a $C^1$ simple closed curve in an oriented surface $M$. Any extremal curve $\alpha \co \R \longrightarrow M \setminus \gamma_0$  positively asymptotic to $\gamma_0$ on the right intersects transversally any extremal curve $\alpha \co \R \longrightarrow M \setminus \gamma_0$  negatively asymptotic to $\gamma_0$ on the right.
\end{lemma}
\proof
Let
\begin{itemize}
	\item $\alpha_0 \co \R \longrightarrow M \setminus \gamma_0$  be a $C^1$ curve positively asymptotic to $\gamma_0$ on the right
	\item $\alpha_1 \co \R \longrightarrow M \setminus \gamma_0$  be a $C^1$ curve negatively asymptotic to $\gamma_0$ on the right
	\item $\delta$ be a $C^1$ transverse segment to $\gamma_0$, oriented so that its transverse intersection with $\gamma_0$ is positive.
\end{itemize}
Intersecting transversally with a given sign is an open property, so there exists a neighborhood $U$ of $\left(\gamma_0, \dot{\gamma_0}\right)$ in $TM$ such that for any $C^1$ arc $\alpha$ in $M$, if $\left(\alpha (t), \dot{\alpha}(t)\right)$ is contained in $U$ for a sufficiently long time, then $\alpha$ intersects $\delta$ transversally with positive  sign.

Since $\alpha_1$ is negatively asymptotic to $\gamma_0$ on the right, there exists a tubular neighborhood $V$ of $\gamma_0$ in $M$ such that $\alpha_1$ eventually leaves the right-hand side of $V$. Restricting $U$ if necessary, we assume $p (U) \subset V$. Take 
 $t_1,t_2$  two consecutive intersection points of $\alpha_0$ with $\delta$, such that 
 $\alpha_0 (\left[t_1,t_2\right]) \subset V$.

Consider the topological annulus $A$ bounded by $\gamma_0$ on the left, and on the right, by $\alpha_0 (\left[t_1,t_2\right])$ glued to the segment of $\delta$ comprised between $\alpha_0 (t_2)$ and  $\alpha_0 (t_1)$. Since $\alpha_1$ eventually leaves the right-hand side of $V$, it must leave $A$. In so doing it cannot intersect $\delta$ for then the intersection of $\delta$ with $\alpha_1$ would be negative. Therefore it must intersect $\alpha_0$, which proves the lemma.
\qed
\subsection{Periodic orbits which are not fixed points}
Besides the Aubry set, another set of note is the Ma\~{n}\'e set $\tilde{\mathcal{N}}(L)$ ; all we need to know about it is that 
\begin{itemize}
	\item it is compact and $\phi_t$-invariant
	\item it contains $\tilde{\mathcal{A}}(L)$
	\item no projection to $M$ of an  orbit contained in $\tilde{\mathcal{N}}(L)$ intersects tranversally the projection to $M$ of an orbit contained in $\tilde{\mathcal{A}}(L)$ (\cite{Fathi}, Theorem 5.2.4)
	\item it is lower semi-continuous with respect to the Lagrangian, that is, for any neighborhood $U$ of $\tilde{\mathcal{N}}(L)$, for any sequence $L_n$ of Tonelli Lagrangians converging to $L$ in the $C^2$ compact-open topology, for $n$ large enough $\tilde{\mathcal{N}}(L_n) \subset U$.
\end{itemize}
\begin{lemma}\label{local}
Let 
\begin{itemize}
\item
$c$ be a cohomology class 
\item
 $\gamma_i, i \in I$ be a collection of $c$-minimizing periodic orbits which are not fixed points
 \item
 $\mu_i$ be the minimizing probability supported on $\gamma_i, i \in I$
 \item $h_i$ be the homology class  of $\mu_i$, for all $i \in I$ 
  \item 
  $\bar{h}$ be any barycenter, with positive coefficients, of the homology classes $h_i$  \item 
  $F :=   \mathcal{L}(\bar{h})$.
 \end{itemize}
 Then there exists  a neighborhood $V$ of $\cup_{i \in I}\gamma_i$ in $TM$, such that $\tilde{\mathcal{A}} (F) \cap V$ consists of closed orbits whose projections to $M$ are homologous to $\gamma_i$ for some $i$.
\end{lemma}
\proof 

 The face $F= \mathcal{L}(\bar{h})$ of $\alpha$  contains $c$ (not necessarily in its interior) because the curves $\gamma_i$ are  $c$-minimizing.  Choose one of the $\gamma_i$; assume, without loss of generality, that it is $\gamma_0$ and its homology class is $h_0$. 

\textbf{Case 1 : $M$ is oriented.}

Denote by $\inter$ the symplectic intersection form induced on $H_1 (M,\R)$ by the algebraic intersection number of closed curves.

\textbf{Case 1.1 : $\gamma_0$ does not separate $M$.}
Before launching into the proof, let us explain the idea. By Lemma \ref{key_o}, there exists a neighborhood $U$ of $\gamma_0$ such that $\mathcal{A}(F) \cap U$ consists of 
periodic orbits homotopic to $\gamma_0$, and of orbits asymptotic to periodic orbits homotopic to $\gamma_0$. What we want is to exclude the asymptotic orbits.

Assume we can find  cohomology classes $c^+, c^-$ in $F$, such that 
\begin{itemize}
\item
$\mathcal{N}(c^+)$ contains orbits positively asymptotic on the right to $\gamma_0$
\item
$\mathcal{N}(c^+)$ contains orbits negatively asymptotic on the left to  $\gamma_0$
\item
$\mathcal{N}(c^-)$ contains orbits positively asymptotic on the left to $\gamma_0$
\item
$\mathcal{N}(c^-)$ contains orbits negatively asymptotic on the right  to $\gamma_0$.
\end{itemize}
Observe that 
$$
\mathcal{A}(F) \subset \mathcal{A}(c^{\pm}) 
$$
since $c^{\pm} \in F$. Therefore no orbit contained in $\mathcal{A}(F)$ can intersect 
$\mathcal{N}(c^+)$ or $\mathcal{N}(c^-)$. On the other hand, any orbit asymptotic to 
$\gamma_0$ must intersect either $\mathcal{N}(c^+)$ or $\mathcal{N}(c^+)$ by Lemma \ref{asymptotic}. Therefore no orbit asymptotic to $\gamma_0$ is contained in $
\mathcal{A}(F)$. 

In real life such $c^{\pm}$ need not exist, because there may be minimizing periodic orbits homotopic to $\gamma_0$ that accumulate on $\gamma_0$. Still, essentially the same idea works. Now let us start the proof.

Since $\gamma_0$ does not separate $M$,  $h_0 \neq 0$, so by the non-degeneracy of $\inter$, we may pick $h' \in \Gamma $ such that $\inter(h_0,h') =1$. Denote $h^{\pm}_{n}:= n\bar{h} \pm h' \in \Gamma , n \in \N $. Take, for each $n \in \N^* $, 
\begin{itemize}
	\item $n\inv h^{\pm}_{n}$-minimizing measures $\mu^{\pm}_{n}$ 
	\item $c^{\pm}_{n} \in H^1 (M,\R)$ such that $\mu^{+}_{n}$ (resp. $\mu^{-}_{n}$) is $c^{+}_{n}$-minimizing (resp. $c^{-}_{n}$-minimizing).
\end{itemize}
The homology classes  $n\inv  h^{\pm}_{n}$ remain in a compact subset  of  $H_1 (M,\R)$ so the cohomology classes $c^{\pm}_{n}$ remain in a compact subset  of  $H^1 (M,\R)$. Therefore  the supports of the measures $\mu^{\pm}_{n}$, which lie in the energy levels $\alpha (c^{\pm}_{n})$ by \cite{Carneiro},  remain in a compact subset  of  $TM$. Hence the sequences $\mu^{\pm}_{n}$, $n\in \N$, have weak$^{\ast}$ limit points $\mu^{\pm}$. Likewise, we may assume the sequences $c^{\pm}_{n}$ converge to $c^{\pm} \in H^1 (M,\R)$. Since the homology class is a continuous function of the measure, we have $\left[\mu^{\pm}\right]=\bar{h} $.  Besides, since $\mu^{\pm}_{n}$ is $c^{\pm}_{n}$-minimizing, 
	\begin{eqnarray*}
\langle c^{+}_{n},	h^{+}_{n}\rangle   &=& \alpha (c^{+}_{n}) + \beta (h^{+}_{n})\\
\langle c^{-}_{n},	h^{-}_{n}\rangle   &=& \alpha (c^{-}_{n}) + \beta (h^{-}_{n})
\end{eqnarray*}
whence, taking limits,
	\[\langle c^{\pm},	h\rangle   = \alpha (c^{\pm}) + \beta (h)
\]
that is, $c^{\pm} \in \mathcal{L}(\bar{h})$.

By Mather's Graph Theorem (\cite{Mather91}, see also \cite{Fathi}, Theorem 5.2.4) there exists $K>0$ such that 
$$
\forall (x,v) \in \tilde{\mathcal{A}}_0,\  \forall (x,v) \in 
\tilde{\mathcal{N}}(c^+), \ 
d_{TM} \left( (x,v),(x',v') \right) \leq K d_M (x,x').
$$
Thus if $V$ is a neighborhod of $(\gamma_0,\dot\gamma_0)$ in $TM$, given by Lemma \ref{key_o}, there exists a neighborhood $U$ of $\gamma_0$ in $M$, such that $U \subset p(V)$, and, for all $x$ in $\overline{U} \cap\mathcal{N}(c^+)$; for all $v \in T_x M $ such that $(x,v) \in \tilde{\mathcal{N}}(c^+)$, we have $(x,v) \in V$.

Let us take $U$ as above, and such that, furthermore, the closure $\overline{U}$ is diffeomorphic to an annulus. Since $\inter (h_0, \left[ \mu_n^+\right])>0$, for any  $n \in \N^* $ there exists an orbit segment $X_n$, contained in $\spt \mu_n^+$, such that $p(X_n)$ crosses $U$ from right to left with respect to the orientation of $U$ induced by that of $M$. 

Denote by $X$ some limit point, with respect to the Hausdorff distance on compact sets, of the sequence $X_n$.
By the semi-continuity of the Ma\~n\'e set with respect to the Hausdorff distance on compact sets, we have
$$ X \subset \tilde{\mathcal{N}}(c^+) \cap V$$
so $X$ consists of periodic orbits homotopic to $\gamma_0$, and of orbits (or pieces of orbits) asymptotic to 
periodic orbits homotopic to $\gamma_0$. Denote by $G$ the set of periodic orbits homotopic to $\gamma_0$ that are contained in $X$.

First observe that any element of $G$ carries an invariant measure, so if $\gamma \in G$, $\gamma (\R) \subset \mathcal{M}(c^+)$. Therefore, by Mather's Graph Theorem, no two elements of $G$ intersect. Thus if $\gamma, \gamma' \in G$ either $\gamma$ lies entirely on the right of $\gamma'$, or $\gamma$ lies entirely on the left of $\gamma'$. We write $\gamma> \gamma'$ if $\gamma$ lies entirely on the left of $\gamma'$. The set $G$ is thus totally ordered by the relation $>$.

Recall that a $(T,\epsilon)$-pseudo orbit is a piecewise continuous curve made up with portions of orbit of $\phi_t$, defined on intervals of length $\geq T$, with finitely many discontinuities, such that at each discontinuity the jump is smaller than $\epsilon$. 
By Lemma 14 of \cite{Bernard-Conley} (which itself elaborates on ideas of \cite{Conley}), the set $X$ enjoys the following property : for all $x,y$ in $X$, either for all $T,\epsilon >0$ there exists a $(T,\epsilon)$-pseudo orbit from $x$ to $y$, or for all $T,\epsilon >0$ there exists a $(T,\epsilon)$-pseudo orbit from $y$ to $x$.

In particular, if $\gamma \in G$ has a successor $\gamma'$ with respect to the order $>$, then there exists an orbit in $X$ positively asymptotic to $\gamma'$  on the right, and negatively asymptotic to $\gamma$  on the left. 

Now let us show that $G$ contains all the  periodic orbits homotopic to $\gamma_0$ and contained in $\mathcal{A}(c) \cap U$. Let $\gamma$ be such a periodic orbit. Since 
$\gamma$ is contained in $\mathcal{A}(c) \subset \mathcal{A}(c^+)$, $\gamma$ does not intersect transversally the projection of any orbit contained in $\mathcal{N}(c^+)$, hence the order $>$ is well-defined on 
$G \cup \{ \gamma\}$. Assume the sets $ \{ \gamma' \in G \co \gamma' < \gamma \}$ and $ \{ \gamma' \in G \co \gamma' > \gamma \}$ are both non-empty.
Let 
\begin{eqnarray*}
\gamma_1 & := & \max \{ \gamma' \in G \co \gamma' < \gamma \}\\
\gamma_2 & := & \min \{ \gamma' \in G \co \gamma' > \gamma \}.
\end{eqnarray*}
The minimum and maximum are well-defined thanks to the compactness of $X$, and $\gamma_2$ is the successor of $\gamma_1$. So there exists a heteroclinic orbit in $X$ from 
$\gamma_1$ to $\gamma_2$. Assume $\gamma \neq \gamma_1$ and $\gamma \neq \gamma_2$, then the heteroclinic orbit from $\gamma_1$ to $\gamma_2$ intersects $\gamma$, which contradicts the fact that $\gamma$ does not intersect 
$\mathcal{N}(c^+)$. Therefore $\gamma = \gamma_1$ or $\gamma = \gamma_2$, in particular $\gamma \in G$, in the case where $ \{ \gamma' \in G \co \gamma' < \gamma \}$ and $ \{ \gamma' \in G \co \gamma' > \gamma \}$ are both non-empty.

Observe that $ \{ \gamma' \in G \co \gamma' < \gamma \}$ and $ \{ \gamma' \in G \co \gamma' > \gamma \}$ cannot both empty. Now assume that $ \{ \gamma' \in G \co \gamma' < \gamma \}$, for instance, is non-empty. Let $\gamma_1$ be its maximum. 
Assume $\gamma \neq \gamma_1$. There exists a heteroclinic orbit positively asymptotic to $\gamma_1$ on the right that comes from the boundary of $U$. This orbit  
must intersect $\gamma$, which is again a contradiction.
 
 So $G$ contains all the  periodic orbits homotopic to $\gamma_0$ and contained in $\mathcal{A}(c) \cap U$.

Now assume that $\mathcal{A}(c) \subset \mathcal{A}(c^+)$ contains an orbit $\delta$ positively asymptotic on the left to a periodic orbit $\gamma$ homotopic to $\gamma_0$.
Then if $\gamma$ has a successor in $G$, $\delta$ intersects the heteroclinic orbit from $\gamma$ to its successor, which contradicts the fact that $\delta$ is contained in $\mathcal{A}(c)$. If $\gamma$ does not have a successor in $G$, then $\delta$ intersects one of the periodic orbits that accumulate on the left of $\gamma$.

In any case we have proved that  $\mathcal{A}(c) \cap U$ does not contain any orbit  positively asymptotic on the left to a periodic orbit $\gamma$ homotopic to $\gamma_0$.
A similar argument proves that $\mathcal{A}(c) \cap U$ does not contain any orbit  negatively asymptotic on the right to a periodic orbit $\gamma$ homotopic to $\gamma_0$. 

Reasoning with $\mu_n^-$ instead of $\mu_n^+$, we show in the same fashion that 
$\mathcal{A}(c) \cap U$ does not contain any orbit  positively asymptotic on the right (resp. negatively asymptotic on the left) to a periodic orbit $\gamma$ homotopic to $\gamma_0$.

We have proved that $\mathcal{A}(c) \cap U$ consists of  periodic orbits homotopic to $\gamma_0$.

\textbf{Case 1.2 : $\gamma_0$  separates $M$.}
\begin{remark}
In that case  the result is stronger : there exists a  neighborhood of $(\gamma_0, \dot{\gamma}_0)$ in $TM$ such that $\tilde{\mathcal{A}}(c) \cap V$ consists of closed orbits whose projection to $M$ are homotopic to $\gamma_0$ (recall that 
$\tilde{\mathcal{A}}\mathcal{L}( \bar{h} )   \subset \tilde{\mathcal{A}}(c)$ because 
$ c\in \mathcal{L}( \bar{h} )$).
\end{remark}

Let $V$  be the neighborhood of $(\gamma_0, \dot{\gamma}_0)$ in $TM$ given by Lemma \ref{key_o}, such that   for any orbit $\phi_t (x,v)$ in $\tilde{\mathcal{A}}(c)$, if 
$\phi_t (x,v)$ meets $V$, then $p\circ \phi_t (x,v)$ is either
\begin{itemize}
  \item   a closed orbit homotopic to $\gamma_0$
    \item or asymptotic to  a closed orbit homotopic to $\gamma_0$.
\end{itemize}
In the former case, we are done, so assume that for some $(x_0,v_0)$ in $\tilde{\mathcal{A}}(c) \cap V$, $\gamma(t) := p\circ \phi_t (x_0,v_0)$ is asymptotic to 
 a closed orbit $\gamma_1$ homotopic to $\gamma_0$.
 
 For definiteness we assume $\gamma(t)$ is positively asymptotic on the left to 
  $\gamma_1$, the three other cases are identical.
  
  Let $\delta \co \left] -1,1\right[ \longrightarrow M$ be a geodesic  arc such that $\delta(0)=\gamma_1(0)$ and $\Omega(\dot\gamma_1 (0), \dot\delta(0)) > 0$, where $\Omega$ is the orientation two-form of $M$. 
  
Then, since $\gamma(t)$ is  positively asymptotic to $\gamma_1$, $\gamma$ intersects $\delta$ infinitely many times, so there exist a sequence $t_n$ of real numbers such that $\lim t_n = +\infty$ and $\gamma(t_n)=\delta(s_n)$ for some $s_n \in  \left] -1,1\right[ $.

By Mather's Graph Theorem, there exists $\epsilon_1 > 0$ such that for any $(x,v)$ in 
$\tilde{\mathcal{A}}(c) $ with $d (x, \gamma_1) < \epsilon_1$ and $x=\delta(s)$ for some $s \in  \left] -1,1\right[ $, we have $\Omega(v, \dot\delta(s)) > 0$. Set
$$ \epsilon_2 := \min \{ \epsilon_1, \frac{1}{2}d(x_0, \gamma_1)\}.$$
Take $n \in \N$ such that $\gamma(\left[t_n,t_{n+1}\right])$ is contained in the 
$\epsilon_2$-neighborhood of $\gamma_1$. Let $D$ be the closed domain bounded on the right by $\gamma_1$, and on the left by the closed curve obtained by joining 
$\gamma(\left[t_n,t_{n+1}\right])$ with the segment of $\delta$ comprised between 
$\gamma(t_n)$ and $\gamma(t_{n+1})$. Observe that $x_0 \not\in D$ because we have required that $\epsilon_2 < d(x_0, \gamma_1)$.

Now pick $(x,v)$ in  $\tilde{\mathcal{A}}(c)$, such that $x \in D$. The curve 
$\alpha(t) := p\circ \phi_t (x,v)$ cannot leave $D$ through $\gamma$ or $\gamma_1$, by the Graph Theorem; nor can it leave through $\delta$, because at any intersection point of $\delta$ and $\alpha$ we must have $\Omega(\dot\alpha, \dot\delta) > 0$. Thus $\alpha(t) \in D$ for all $t\geq 0$.

Furthermore $(x,v) \in V$ so $\alpha$ is either
\begin{itemize}
  \item   a closed orbit homotopic to $\gamma_0$
    \item or asymptotic to  a closed orbit homotopic to $\gamma_0$.
\end{itemize}
Note that if $\alpha$ were closed and distinct from $\gamma_1$, then $\alpha$ would lie in $D$ on the left of $\gamma_1$, so $\gamma$ would intersect $\alpha$. Likewise if $\alpha$ were asymptotic to a  closed orbit  $\alpha_1$ distinct from $\gamma_1$, then $\gamma$ would intersect $\alpha_1$. In any case this would be a contradiction to  the Graph Theorem, so we have proved that for any $(x,v)$ in  $\tilde{\mathcal{A}}(c)$, such that $x \in D$, the curve 
$\alpha(t) = p\circ \phi_t (x,v)$ is either $\gamma_1$ itself or is positively asymptotic to $\gamma_1$. 

Denote by $M_2$ the connected component of $M\setminus \gamma_1$ that does not contain $x_0$.
Take $\epsilon > 0$ such that the $2\epsilon$-neighborhood of $\gamma_1$ is contained in $D \cup M_2$. 

Since $D$ is compact, there exists $T_1 \in \R$ such that 
for any $(x,v)$ in  $\tilde{\mathcal{A}}(c)$, with $x \in D$, $\forall t \geq T_1$, still denoting $\alpha(t) = p\circ \phi_t (x,v)$, we have $d(\alpha(t), \gamma_1) < \epsilon$. 

Take $T_2$ such that $d(\gamma(t), \gamma_1)<\epsilon$ for all $t \geq T_2$, and set 
$$T := \max \{T_1, T_2\}.$$

Recall from \cite{Fathi} that the Aubry set is chain-recurrent, so there exists  finite sequences $(x_1,v_1), \ldots (x_n,v_n)=(x_0,v_0)$ and $t_1, \ldots t_n$ such that for $i=1,\ldots n$
\begin{itemize}
  \item $(x_i,v_i) \in \tilde{\mathcal{A}}(c)$ 
   \item $t_i > T$ 
    \item $d_{TM} \left( \phi_ {t_i} (x_{i-1}, v_{i-1}), (x_i,v_i)\right) < \frac{\epsilon}{K}$, where $K$ is given by the Graph Theorem.
\end{itemize}
Since $t_1 > T \geq T_2$, we have $d(\gamma(t), \gamma_1) < \epsilon$. Moreover,
$$d(x_1, \gamma(t_1) ) \leq d_{TM}  \left( \phi_ {t_i} (x_{0}, v_{0}), (x_1,v_1)\right) < \epsilon
$$
so  $d(x_1, \gamma_1)< 2\epsilon$, whence $x_1 \in D \cup M_2$. Then\begin{itemize}
  \item if $x \in D$, we have 
  $$d( p \circ \phi_{t_2}(x_1,v_1), \gamma_1) < \epsilon$$
  because $t_2 > T_1$
  \item 
if $x_1 \in M_2$, then $ \phi_{t_2}(x_1,v_1) \in M_2$, because the curve $ \phi_{t}(x_1,v_1)$, which cannot intersect $\gamma_1$, cannot leave $M_2$ (this is where we use the fact that $\gamma_0$, hence $\gamma_1$, is separating).
\end{itemize}
Thus either $x_2 \in M_2$, or $ d(x_2, \gamma_1) < 2\epsilon$; in any case 
$x_2 \in D\cup M_2$. By induction we prove that $x_n   \in D\cup M_2$. But $x_n = x_0 \not\in D\cup M_2$. 

This contradiction shows there is no $(x_0,v_0)$ in $\tilde{\mathcal{A}}(c) \cap V$, such that $\gamma(t) = p\circ \phi_t (x_0,v_0)$ is asymptotic to 
 a closed orbit  homotopic to $\gamma_0$.

Thus  $\tilde{\mathcal{A}}(c) \cap V$ consists of closed orbits whose projection to $M$ are homotopic to $\gamma_0$.

Now let us finish the proof of the orientable case of the lemma. For every $i \in I$, we have found a neighborhood $V_i$ of $(\gamma_i,\dot{\gamma}_i)$ in $TM$ such that $\tilde{\mathcal{A}} \left(\mathcal{L}(\bar{h} )\right)\cap V_i$ consists of periodic orbits homologous to $\gamma_i$. Define  $V$ to be the union over $i \in I$ of the $V_i$, then $V$ is a neighborhood of $\cup_{i \in I}\left(\gamma_i, \dot{\gamma}_i\right)$ in $TM$, and  $\tilde{\mathcal{A}} \left(\mathcal{L}(\bar{h} )\right)\cap V$ consists of periodic orbits homologous to $\gamma_i$ for some $i$.

\textbf{Case 2 : $M$ is not orientable.}
Let $\pi \co M_o \longrightarrow M$ be the orientation cover of
$M$. Then $M_o$ is an orientable surface endowed with a
fixed-point free, orientation-reversing involution $I$. As in Appendix B, we call $E_1$ (resp $E_{-1}$) the eigenspace of the homomorphisms of $H_1 (M_o, \R)$ and $H^1 (M_o, \R)$ induced by $I$. So there is one $E_1$ in $H_1 (M_o, \R)$ and another in $H^1 (M_o, \R)$; it will always be clear from the context which is which. Let
\begin{itemize}
	\item $\delta_j$, $j \in J$ be the collection of all lifts to $M_o$ of the $\gamma_i$
	\item $\left\{h'_1,\ldots h'_l\right\}$ be the set of all homology classes of the $\delta_j$, $j \in J$
	\item $\bar{h}'$  be $ (h'_1 +\ldots + h'_l)/l$
	\item $\bar{h}$  be $\pi_* (\bar{h}')$.
\end{itemize}
The set  $\left\{\delta_j \co j \in J \right\}$ is invariant under $I$, thus the set $\left\{h'_1,\ldots h'_l\right\}$ is invariant under $I_*$. Therefore $\bar{h}' \in E_1 \subset H_1(M_o,\R)$. The $\delta_j$, $j \in J$ are minimizing (\cite{Fathi98}) so by the orientable case of the lemma, there exists a neighborhood $V$ of $\cup_{j \in J}(\delta_j, \dot{\delta}_j )$ in $TM_o$ such that $\tilde{\mathcal{A}} \left(\mathcal{L}(\bar{h}' )\right)\cap V$ consists of periodic orbits homologous to $\delta_j$ for some $j$. Taking a smaller $V$ if we have to, we assume that $V$ is invariant under $I$. Lemma \ref{mathcalLpi_2} says
	\[ \pi^* \left(\mathcal{L}(\bar{h})\right)= \mathcal{L}(\bar{h}') \cap E_1.
\]
Now take $c$ in the relative interior of $\mathcal{L}(\bar{h})$. Then, $\pi*$ being a linear isomorphism onto $E_1 \subset H^1(M,\R)$, $\pi^* c$ lies in the relative interior of  $\mathcal{L}(\bar{h}') \cap E_1$. Since $\bar{h}' \in E_1 \subset H_1(M_o,\R)$, it is easy to see that $\mathcal{L}(\bar{h}')$ is $I^*$-invariant. Therefore $\pi^* c$ lies in the relative interior of  $\mathcal{L}(\bar{h}')$, so 
	\[\tilde{\mathcal{A}} \left(\mathcal{L}(\bar{h}' )\right)=\tilde{\mathcal{A}} (\pi^* c).
\]
By \cite{Fathi98}, denoting by $T\pi$ the tangent map to $\pi$, we have
	\[T\pi \left(\tilde{\mathcal{A}} (\pi^* c)\right)= \tilde{\mathcal{A}} (c)
\]
that is
	\[T\pi \left(\tilde{\mathcal{A}} \left(\mathcal{L}(\bar{h}' )\right)\right)=\tilde{\mathcal{A}} \left(\mathcal{L}(\bar{h} )\right).
\]
Now  observe that 
\[T\pi \left(\tilde{\mathcal{A}} \left(\mathcal{L}(\bar{h}' )\right) \cap V\right) \subset T\pi \left(\tilde{\mathcal{A}} \left(\mathcal{L}(\bar{h}' )\right)\right)\cap T\pi (V).
\]
Then take a fundamental domain $M_1$ for $I$, that is, a subset $M_1$ of $M_o$ such that $\pi$ restricted to  $M_1$ is one-to one and  onto $M$. Then 
\begin{eqnarray*}
T\pi \left(\tilde{\mathcal{A}} \left(\mathcal{L}(\bar{h}' )\right) \cap V \cap TM_1 \right)&=&
T\pi \left(\tilde{\mathcal{A}} \left(\mathcal{L}(\bar{h}' \right) \right)
 \cap T\pi \left( V \right) \cap T\pi (TM_1)\\
 &=& T\pi \left(\tilde{\mathcal{A}} \left( \mathcal{L}(\bar{h}' \right) \right)
 \cap T\pi \left( V \right) \cap TM \\
 &=& T\pi \left(\tilde{\mathcal{A}} \left( \mathcal{L}(\bar{h}' \right) \right)
 \cap T\pi \left( V \right)
\end{eqnarray*}
so 

	\[T\pi \left(\tilde{\mathcal{A}} \left(\mathcal{L}(\bar{h}' )\right) \cap V\right)= T\pi \left(\tilde{\mathcal{A}} \left(\mathcal{L}(\bar{h}' )\right)\right)\cap T\pi (V).
\]
The projection $\pi$ is a local diffeomorphism so $T\pi (V)$ is a neighborhood of $\cup_{i \in I}\left(\gamma_i, \dot{\gamma}_i\right)$ in $TM$. This finishes the proof of Lemma \ref{local}.
\qed

The following corollary of Lemma \ref{local} reduces the proof of our main result to proving that $\mathcal{L}(h)$ is a rational flat when $h$ is 1-irrational and nonsingular.
\begin{corollary}\label{corollaire_local}
Assume that for some $c$ in $H^1 (M,\R)$, the Mather set $\tilde{\mathcal{M}}(c)$ consists of periodic orbits which are not fixed points $\gamma_i, i \in I$. Let $h$ be any barycenter with positive coefficients of the homology classes of $\gamma_i, i \in I$. Then 
	\[\tilde{\mathcal{A}}\left(\mathcal{L}(h)\right)= \tilde{\mathcal{M}}\left(\mathcal{L}(h)\right)= \tilde{\mathcal{M}}(c).
\]
\end{corollary}
\proof
By Lemma \ref{local} there exists a neighborhood $V$ of $\tilde{\mathcal{M}}(c)$, such that 
\begin{equation}\label{cor_loc}
\tilde{\mathcal{A}}\left(\mathcal{L}(h)\right) \cap V = \tilde{\mathcal{M}}(c).
\end{equation}
Hence
	\[\tilde{\mathcal{M}}(c) \subset \tilde{\mathcal{A}}\left(\mathcal{L}(h)\right), \mbox{ so } \tilde{\mathcal{M}}(c) \subset \tilde{\mathcal{M}}\left(\mathcal{L}(h)\right).
\]
On the other hand $c \in \mathcal{L}(h)$, thus 
	\[\tilde{\mathcal{M}}(c) \supset \tilde{\mathcal{M}}\left(\mathcal{L}(h)\right), \mbox{ so } \tilde{\mathcal{M}}(c) = \tilde{\mathcal{M}}\left(\mathcal{L}(h)\right).
\]
Now $\tilde{\mathcal{A}}(\mathcal{L}(h))$ consists of $\tilde{\mathcal{M}}(\mathcal{L}(h))$, and orbits homoclinic to $\tilde{\mathcal{M}}(\mathcal{L}(h))$. Such orbits enter any neighborhood of $\tilde{\mathcal{M}}(\mathcal{L}(h))$, so (\ref{cor_loc}) implies 
	\[\tilde{\mathcal{A}}\left(\mathcal{L}(h)\right)= \tilde{\mathcal{M}}\left(\mathcal{L}(h)\right)= \tilde{\mathcal{M}}(c).
\]

\qed

\section{Proof of the main theorem}
 \subsection{Notations}
 Let $M$ be a closed oriented surface and let $L$ be a Tonelli Lagrangian on $M$.
For $h \in H_1(M,\R)\setminus \left\{0\right\}$, we define the maximal radial flat $R_h$ of $\beta$ containing $h$ as the largest subset of the half-line $\left\{th \co t \in \left[0,+\infty \right[ \right\}$ containing $h$ (not necessarily in its relative interior) in restriction to which $\beta$ is affine. Beware that $R_h$ is  a flat of $\beta$, but may not be a face of $\beta$. The possibility of radial flats is the most obvious difference between the $\beta$ functions of Riemannian metrics (\cite{gafa}, \cite{nonor}) and those of general Lagrangians. An instance of radial flat is found in \cite{Carneiro_Lopez}. We define the Mather set $\tilde{\mathcal{M}}(R_h)$  as  the closure in $TM$ of the union of the supports of all $th$-minimizing measures, for $th \in R_h$.

Let $h$ be a homology class. Assume $h$ is  1-irrational.
Then for all $t$ such that $th \in R_h$, $th$ is also 1-irrational. Furthermore $R_h$ is contained in a face of $\beta$, so Mather's Graph Theorem and Proposition \ref{rational} combine to say that $\tilde{\mathcal{M}}(R_h)$ is a union of pairwise disjoint periodic orbits $\gamma_i$, $i \in I$ where $I$ is some set, not necessarily finite. We denote by $V(R_h)$ the linear subspace of $H_1(M,\R)$ generated by $\left[\gamma_i\right]$, $i \in I$.

 Recall that $\mathcal{L}(h)$ is the Legendre transform of $h$ with respect to $\beta$. In Appendix A we define  $\mathcal{F}(\mathcal{L}(h))$ to be  the set of homology classes $h'$ that lie in $\mathcal{L}(c)$ for all $c \in \mathcal{L}(h)$. By Lemma \ref{interieur_face}, $\mathcal{F}(\mathcal{L}(h))$ is a face of $\beta$. It is clear that $h$ lies in $\mathcal{F}(\mathcal{L}(h))$, although possibly on the boundary.

\subsection{}
\begin{lemma}\label{nonsingular_1}
Let $L$ be a Tonelli Lagrangian on a closed surface $M$.
Let $h$ be a $1$-irrational, nonsingular  homology class. Then the relative interior of $R_h$ is contained in the relative interior of the face  $\mathcal{F}(\mathcal{L}(h))$.
\end{lemma}
\proof
\textbf{Orientable case}

Assume $M$ is oriented. Let us denote $V_0 := V(R_{h})$ for short.
 Take $h_i, i=1\ldots k$ such that $ h +h_i$ lies in $\mathcal{F}(\mathcal{L}(h))$ for all $i=1,\ldots k$, and the convex hull of $ h +h_i$, $i=1,\ldots k$, contains an open subset of $\mathcal{F}(\mathcal{L}(h))$. Pick one of the $h_i$. Then the segment $\left[h,h+h_i\right]$ is contained in a face of $\beta$ so by Mather's graph Theorem $h_i \in V^{\perp}_{0}$.
 
Then $-h_i$ lies in $ V^{\perp}_{0}$ so by Theorem \ref{thm_gafa} there exist $t_i \in \R$, $s_i := s(h,h_i) >0$ such that the segment $\left[ h,t_i h - s_i h_i\right]$ is contained in a face of $\beta$. So there exists $c \in H^1(M,\R)$ such that  
\begin{eqnarray*}
\alpha (c) + \beta ( h) &=& \langle c,h\rangle   \\	
\alpha (c) + \beta (t_i h - s_i h_i) &=& \langle c,t_i h - s_i h_i\rangle  .
\end{eqnarray*}
The first equality above says that $c \in \mathcal{L}(h)$. Since $ h +h_i$ lies in $\mathcal{F}(\mathcal{L}(h))$, we also have 
	\[\alpha (c) + \beta ( h +  h_i) = \langle c, h +  h_i\rangle  .
\]
Therefore 
$\mathcal{F}(\mathcal{L}(h))$ contains the convex hull $C_i$ of $h$, $ h +h_i$, $t_i h - s_i h_i$ for $i=1,\ldots k$.  

We claim that for each $i=1,\ldots k$, $C_i$ contains some $t h$ in its interior. Indeed, the following three cases may occur : 
\begin{itemize}
	\item $h_i \notin \R h,\  t_i \neq 1$ : then $C_i$ is a 2-simplex containing $th$ in its relative interior, for some $t \in \left] 1,t_i \right[$, for instance $t= (1+2s_i +t_i)/(2+2s_i)$
	\item $h_i \notin \R h,\  t_i = 1$ : then $C_i$ is the segment $\left[h +h_i, h - s_i h_i \right]$, which contains $h$ in its relative interior
	\item $h_i \in \R h$ : then $C_i$ is a segment of the straight line $\R h$, hence it contains some $th$ in its relative interior.
\end{itemize}
Now the convex hull $C$ of $\cup^{k}_{i=1} C_i$ is contained in $\mathcal{F}(\mathcal{L}(h))$ by convexity of the latter. The relative interior of $C$ is open in $\mathcal{F}(\mathcal{L}(h))$ by our hypothesis on $h_i, i=1\ldots k$. On the other hand $C$ contains an interior point of $R_h$ in its relative interior, because $C$ is the convex hull of finitely many convex sets, each of which contains an interior point of $R_h$ in its relative interior.

So the intersection of the relative interiors of $R_h$ and $\mathcal{F}(\mathcal{L}(h))$ is non-empty. Then by \cite{Rockafellar}, Theorem 6.5, the intersection of the relative interiors of $R_h$ and $\mathcal{F}(\mathcal{L}(h))$ is the relative interior of  $R_h \cap\mathcal{F}(\mathcal{L}(h)) = R_h$.  

\textbf{Non-orientable case}

Take
\begin{itemize}
	\item $h$  a $1$-irrational, nonsingular  homology class of $M$
	\item $h'\in E_1 \subset H_1 (M_o,\R)$ such that $\pi_* h' =h$
	\item $c \in H^1(M,\R)$ such that $\alpha (c)+\beta (h)= \langle c,h\rangle  $.
\end{itemize}
 Since $\pi_*$ sends an integer class to an integer class, and is one-to-one on $E_1$, $h'$ is $1$-irrational. 

Any support of an $h'$-minimizing measure is the lift to $M_o$ of the support of an $h$-minimizing measure by \cite{Fathi98}, hence $h'$ is nonsingular.

The orientable case of the lemma then says that for any  $t$ such that $th'$ is in the relative interior of $R_{h'}$, $th'$ lies in the relative interior of $\mathcal{F}(\mathcal{L}(h'))$. Furthermore $h' \in E_1$ so $th'$ lies in the relative interior of $\mathcal{F}(\mathcal{L}(h')) \cap E_1$. Thus,  $\pi_*$ being  linear,  $\pi_* (th') =th$ lies in the relative interior of $\pi_{\ast}\left(\mathcal{F}(\mathcal{L}(h'))\cap E_1 \right)$. This combines with Lemma \ref{mathcalLpi} to end the proof.
\qed

A consequence of our last lemma is that for any $1$-irrational, nonsingular  homology class $h$, $\mathcal{L}(h)$ is a nonsingular rational flat of $\alpha_L$ :
\begin{lemma}\label{nonsingular_3}
Let $h$ be a $1$-irrational, nonsingular  homology class. We have 
	\[\tilde{\mathcal{M}}(\mathcal{L}(h))= \tilde{\mathcal{M}}(R_h).
\]
\end{lemma}
\proof
Recall that any measure supported in $\tilde{\mathcal{M}}(\mathcal{L}(h))$is $c$-minimizing, for any $c \in \mathcal{L}(h)$ (see \cite{Mane92}). 
Let $\mu$ be a minimizing measure supported in $\tilde{\mathcal{M}}(\mathcal{L}(h))$. Then for all  $c \in \mathcal{L}(h)$ we have
	\[\alpha (c) + \beta (\left[\mu\right]) = \langle c,\left[\mu\right]\rangle  
\]
whence $\left[\mu\right] \in \mathcal{F}(\mathcal{L}(h))$. Take $t$ such that $th$ lies in the interior of $\mathcal{F}(\mathcal{L}(h))$. Take $\lambda \in \left]0,1\right[$ and $h'$ in $\mathcal{F}(\mathcal{L}(h))$ such that 
	\[ th = \lambda \left[\mu\right] +(1-\lambda) h'.
\]
Take an $h'$-minimizing measure $\mu'$. Then $\lambda \mu +(1-\lambda) \mu'$ is a $th$-minimizing measure, hence its support is contained in $\tilde{\mathcal{M}}(R_h)$.
\qed

\subsection{Proof of Theorem \ref{principal}}
Take a nonsingular, 1-irrational homology class $h$.
Note that $th$ is 1-irrational for any $t$.
 For  any measure $\mu$ supported in $\tilde{\mathcal{M}}(\mathcal{L}(h))$, the homology class of  $\mu$ lies in $\mathcal{F}(\mathcal{L}(h))$. Thus by Lemma \ref{nonsingular_1} and   Lemma \ref{face_rationnelle},   the support of $\mu$ consists of periodic orbits and fixed points. The latter are ruled out by the nonsingularity of $th$, which itself is a consequence of Lemma \ref{R_h_nonsing}.  Thus $\tilde{\mathcal{M}}(\mathcal{L}(h))$ is a union of non-trivial minimizing periodic orbits. Take $c$ in the relative interior of $\mathcal{L}(h)$. We have, by Corollary \ref{corollaire_local},  
	\[\tilde{\mathcal{M}}(\mathcal{L}(h))=\tilde{\mathcal{M}}(\mathcal{L}(th))=\tilde{\mathcal{A}}(\mathcal{L}(th))=\tilde{\mathcal{A}}(\mathcal{L}(h))
\]
which finishes the proof of Theorem \ref{principal}.
\qed

\section{Examples}
\subsection{}
Here is an example to illustrate why we need the non-singularity hypothesis in our theorem. Let $X$ be a vector field on the standard two-sphere $\Sbar^2$ such that\begin{itemize}
  \item $X$ has only one zero at some point $x_0$
  \item $X$ has no periodic orbit
  \item every orbit is asymptotic, positively and negatively, to $x_0$.
\end{itemize}
 Note that every point of $\Sbar^2$ is chain-recurrent under the flow of $X$. Using an idea of Ma\~n\'e let us define a Lagrangian on $T\Sbar^2$ by 
$$
L(x,v):= \frac{1}{2} \| v-X(x)\|^2.
$$
The projected Mather set of $L$ is $ \{x_0\}$, because the only invariant measure of $X$ is the Dirac at $x_0$. 
On the other hand, by \cite{FFR}, Thm 1.6, the projected Aubry set of $L$ is the chain-recurrent set of $X$, which is $\Sbar^2$. To simplify the notation we denote
\begin{itemize}
  \item by $0^1$ the zero element of $H^1 (\Sbar^2,\R)$
  \item  by $0_1$ the zero element of $H_1 (\Sbar^2,\R)$. 
\end{itemize}
Since $H^1(\Sbar^2, \R)= \{0^1\} $, we have 
$$
\mathcal{L}(0_1)=  \{0^1\}
$$
and 
$$
\mathcal{A}(\mathcal{L} (0_1)) =  \Sbar^2 \neq  \{x_0\}=  \mathcal{M}(\mathcal{L} (0_1)).
$$ 
\subsection{}
The last example is a bit too easy, because the sphere has no homology, so here is an example on the two-torus. Let $\T^2$ be  the standard two-torus $\T^2$ and let $(e_1,e_2)$ be a basis of $H_1 (\T^2, \Z)$. We again use  the  notation of the previous paragraph :
\begin{itemize}
  \item  $0^1$ is the zero element of $H^1 (\T^2,\R)$
  \item   $0_1$ is the zero element of $H_1 (\T^2,\R)$. 
\end{itemize}
Let $X$ be a vector field on $\T^2$  such that $X$  vanishes at some point $x_0$, every orbit  of $X$ is homoclinic to $x_0$,  and $X$ has four homoclinic orbits 
$\alpha_1,\alpha_2, \alpha_3,\alpha_4$ to $x_0$ with the following property : denote 
$\alpha'_i$, $i=1,2,3,4$, the closed curve defined on $\left[-\pi /2, \pi/2 \right]$ such that 
$\alpha'_i (t)= \alpha (\tan (t))$ for $t\in \left]-\pi /2, \pi/2 \right[ $, and $\alpha'_i (\pm \pi /2))= x_0$. We require that
\begin{itemize}
  \item  $\alpha'_1$ be homologous to $e_1$
  \item  $\alpha'_2$ be homologous to $-e_1$  
  \item $\alpha'_3$ be homologous to $e_2$
   \item $\alpha'_4$ be homologous to $-e_2$.  
   \end{itemize}
Let  $L$ be the Lagrangian on $T\T^2$ defined by by 
$$
L(x,v):= \frac{1}{2} \| v-X(x)\|^2.
$$
As in the previous example, we have $\alpha_1,\alpha_2, \alpha_3,\alpha_4 \subset \mathcal{A}(L)= \mathcal{A}(0^1)$, but here  since the torus has non-trivial homology, 
$\mathcal{A}(\mathcal{L} (0_1))$ could be properly contained in $\mathcal{A}(0^1)$.
However, we shall prove that $\mathcal{L} (0_1)= \{0^1\}$, so $\mathcal{A}(\mathcal{L} (0_1))= \mathcal{A}(0^1)$. Therefore 
$$
\alpha_1,\alpha_2, \alpha_3,\alpha_4 \subset \mathcal{A}(\mathcal{L} (0_1)),
$$
on the other hand 
$$
\alpha_1,\alpha_2, \alpha_3,\alpha_4 \not\subset \mathcal{M}(\mathcal{L} (0_1))
$$
therefore
$$
\mathcal{M}(\mathcal{L} (0_1)) \neq \mathcal{A}(\mathcal{L} (0_1)).
$$
Now let us prove that $\mathcal{L} (0_1)= \{0^1\}$. Recall that  $\mathcal{L} (0_1)\ni \{0^1\}$ because the Dirac measure  at $(x_0, 0)$, whose homology class is $0_1$, is $0^1$-minimizing. So   proving that $\mathcal{L} (0_1)= \{0^1\}$ amounts to proving that $\mathcal{L} (0_1)$ has only one element, in other words, that  $\beta := \beta_L$ is differentiable at $0_1$.
   
For $n$ large enough there exists a unique geodesic segment $c_n$, parametrizeed with unit speed, between $\alpha (-n)$ and $\alpha (n)$, because both  $\alpha (-n)$ and $\alpha (n)$ converge to $x_0$.  Let $d_n$ be the distance between $\alpha (-n)$ and $\alpha (n)$. Consider the piecewise $C^1$ closed curve  
$$
\begin{array}{rcl }
 \alpha_{1,n} \co \left[0,2n+d_n \right]     &   \longrightarrow  & \T^2\\
   t   &   \longmapsto & \alpha_1 (t-n) \mbox{ for } t \in \left[0,2n\right]\\
    t   &   \longmapsto & c(t-2n) \mbox{ for } t \in \left[2n, 2n+d_n\right].
    \end{array}
$$
For $n$ large enough the homology class of $\alpha_{1,n}$ is $e_1$. Recall from \cite{Bangert99} that a probability measure $\mu$ on $T\T^2$ is said to be closed if $\int df d\mu = 0$ for any $C^1$ function on $M$. The probability measure $\mu_n$ uniformly distributed on  $ \alpha_{1,n}$ is closed because  $ \alpha_{1,n}$ is a  piecewise $C^1$ closed curve. The homology  class of $\mu_n$
is $(2n+d_n)\inv e_1$ (for $n$ large enough).

Since $L$ vanishes on $\alpha_1$, the action of $L$ on  $\mu_n$ is 
$$
\int L d\mu_n = \frac{1}{2n+d_n}\int_0 ^{d_n} \frac{1}{2} \| \dot c (t) -X(c(t))\|^2 dt 
\leq \frac{K d_n}{2n+d_n},
$$
where 
$$
K :=\left(1+ \max_{x \in \T^2 } \|X(x)\|\right)^2 .
$$
Since $\mu_n$  is closed, by \cite {FS}, Theorem 1.6, we have $\beta (\left[\mu_n\right]) \leq \int L d\mu_n$, so, by convexity of $\beta$, $\beta (te_1)= o(t)$ for $t\geq 0$. 
Applying the same method with $\alpha_2$ (resp. $\alpha_3$, resp. $\alpha_4$) instead of $\alpha_1$, we obtain  $\beta (te_1)= o(t)$ for $t\leq 0$ (resp.   $\beta (te_2)= o(t)$ for $t\geq 0$, resp. $\beta (te_2)= o(t)$ for $t\leq 0$). By the convexity of $\beta$, this proves that  $\beta$ is differentiable at zero.
\qed

\appendix

\section{Convex and superlinear functions}
Let
\begin{itemize}
	\item $E$ be a finite dimensional Banach space
	\item $A \co E \longrightarrow \R$ be a convex and superlinear map.
\end{itemize}
Then the Fenchel transform of $A$, defined by the formula
	\[ \begin{array}{rcl}
	B \co E^* & \longrightarrow & \R \\
	y & \longmapsto & \sup_{x \in E}\left(\langle y,x\rangle  -A(x)\right)
	\end{array}
\]
is well-defined, convex and superlinear. The Legendre transform (with respect to $A$) of a point $x$ in $E$ is the set $\mathcal{L}(x)$ of $y \in E^*$ which achieve the supremum above.   Since $B$ is convex and superlinear, there is a Legendre transform with respect to $B$ as well.  We call 
\begin{itemize}
    \item relative interior of a convex subset $C$ of $E$, the interior of $C$ in the affine subspace of $E$ generated by $C$ (see \cite{Rockafellar}, p.44).
    For instance, the relative interior of the interval $\left[a,b\right]$ is $\left\{a\right\}$ if $a=b$, $\left]a,b\right[$ otherwise
    \item supporting subspace to the graph of $A$ any affine subspace of $E \times \R$ that meets the graph of $A$  but not the open epigraph 
    \[\left\{(x,t)\in E\times \R \co t> A(x)       \right\}\]
    \item flat of   $A$, the projection to $E$ of the intersection of the graph of $A$ with a supporting subspace
    \item dimension of a flat, the dimension of the affine subspace it generates in $E$
    \item interior of a flat, its relative interior as a convex set
    \item face of  the graph of $A$, the projection to $E$ of the intersection of the graph of $A$ with a supporting \textit{hyperplane}.
\end{itemize}
Note that 
\begin{itemize}
  \item by the Hahn-Banach Theorem, any flat is contained in a face
  \item for any $x \in E$ (resp. $x \in E^*$), the Legendre transform $\mathcal{L}(x)$ is a face of the graph of $B$ (resp. $A$).
\end{itemize}

Conversely,  by the Hahn-Banach Theorem,  a face of $A$ (resp. $B$) is  the Legendre transform, with respect to $B$ (resp. $A$), of a point of $E^*$ (resp. $E$), that is, a subset $F$ such that
	\[
\exists y_0 \in E^*,\  	F = \left\{ x \in E \co A(x)+B(y_0) =\langle y_0,x\rangle   \right\}.
\]
Given a flat $F$ of $A$, we denote
	\[ \mathcal{F}(F) := \left\{y \in E^* \co \forall x \in F,\  A(x)+B(y) =\langle y,x\rangle  \right\},
\]
that is, $\mathcal{F}(F)$ is the intersection of all  Legendre transforms of points of $F$. Note that for any two flats $F_1,F_2$ of $A$, $F_1 \subset F_2$ is equivalent to $\mathcal{F}(F_2)\subset \mathcal{F}(F_1)$. In particular, if $x,y$ are points of $E^*$, 
\[ \left(x \in \mathcal{F}(\mathcal{L}(y))\right)\Longleftrightarrow  
\left( \mathcal{L}(y) \subset \mathcal{L}(x) \right)\Longleftrightarrow 
\left(\mathcal{F}(\mathcal{L}(x))) \subset\mathcal{F}(\mathcal{L}(y))\right).
\]

\begin{lemma}\label{interieur_face}
Let $F$ be a flat of $A$ and let $x_0$ be a point in the relative interior of $F$. Then $\mathcal{F}(F)$ is the Legendre transform of $x_0$. In particular $\mathcal{F}(F)$ is a face of $B$.
\end{lemma}
\proof
By definition of $\mathcal{F}(F)$, it is contained in  the Legendre transform of $x_0$. Let us show the converse inclusion holds true. Take $y$ such that $A(x_0)+B(y) =\langle y,x_0\rangle  $. We want to show that $y \in \mathcal{F}(F)$, that is, $A(x)+B(y) =\langle y,x\rangle  $ for all $x \in F$. Take $x \in F$. Since $x_0$ lies  in the interior of $F$, there exists $x'$ in $F$ and $0<t<1$ such that $x_0=tx+(1-t)x'$. 
Since any flat is contained in a face, there exists $y_0$ such that $F \subset \left\{ x \in E \co A(x)+B(y_0) =\langle y_0,x\rangle   \right\}$, so 
\begin{eqnarray*}\label{a}
A(x)+B(y_0) &=&\langle y_0,x\rangle  \\
A(x')+B(y_0) &=&\langle y_0,x'\rangle  	.
\end{eqnarray*}
Summing $t$ times the first equation with $(1-t)$ times the second equation yields $tA(x)+(1-t)A(x')+B(y_0) =\langle y_0,x_0\rangle  $, but since $x_0 \in F$, we have $A(x_0)+B(y_0) =\langle y_0,x_0\rangle  $ whence $A(x_0)=tA(x)+(1-t)A(x')$.

On the other hand by definition of $B$ we have
\begin{eqnarray*}\label{b}
A(x)+B(y) &\geq&\langle y,x\rangle  \\
A(x')+B(y) &\geq&\langle y,x'\rangle  	.
\end{eqnarray*}
Summing  $t$ times the first inequality with $(1-t)$ times the second inequality yields the equality $A(x_0)+B(y) =\langle y,x_0\rangle  $, thus both inequalities are equalities, which proves the lemma.
\qed

Our next lemmas are improvements of Lemmas 4.1  and 4.2 of \cite{nonor}.
\begin{lemma}\label{extension_plat}
Let 
\begin{itemize}
  \item $x_0$ be a point of $E$
  \item $I$ be some (possibly infinite) set
  \item $F_i$, $i\in I$ ba a family of flats of $A$, such that
  \item $x_0 \in F_i$ for all $i\in I$
  \item there is at most one index $i$ such that $x_0$ does not lie in the relative interior of $F_i$.
  \end{itemize}  
  Then there exists a flat $F$ containing $F_i$ for all $i\in I$.
\end{lemma}
\proof
Assume for convenience that the one flat that need not contain $x_0$ in its interior is $F_0$.
Let $y \in E^*$ be  such that for all $x$ in $F_0$, we have $B(y)+A(x)=\langle y,x\rangle  $. In particular we have $B(y)+A(x_0)=\langle y,x_0\rangle  $ so $y$ lies in the Legendre transform of $x_0$. Thus by Lemma \ref{interieur_face} $y \in \mathcal{F}(F_i)$ for all $i \in I$. This means that for all $x$ in $F_i$, we have $B(y)+A(x)=\langle y,x\rangle  $. Thus the face $\mathcal{L}(y)$ contains  $F_i$ for all $i\in I$.
\qed

\begin{lemma}\label{nonor_4.2}
Let 
\begin{itemize}
  \item $x_0$ be a point of $E$
  \item $I$ be some (possibly infinite) set
  \item $F_i$, $i\in I$ be a family of flats of $A$ such that $x_0$ lies in the relative interior of $F_i$ for all  $i \in I$.
  \end{itemize} 
  Then there exists a flat $F$ containing $F_ i$ for all $i \in I$ such that  $x$ is an interior point of $F$.
\end{lemma}
By Lemma \ref{extension_plat}, there exists a flat $F_1$ containing $F_i$ for all $i \in I$. In particular $F_1$ contains the convex hull $C_1$ of the union of $F_i$ over all $i\in I$. So the affine subset $V_I$ generated by $C_1$ is a supporting subspace to the graph of $A$. In particular, the intersection of $V_I$ with the graph of $A$ is a flat $F$ of $A$, and $F_i \subset F$ for all $i \in I$. 

Now, since the dimension of $E$ is finite, there exists a finite subset $J$ of $I$ such that $V_I$ is generated by the union of $F_j$, $j \in J$. Since $x_0$ lies in the relative interior of $F_j$ for all $j \in J$, $x_0$ lies in the relative interior of the convex hull  $C_2$ of the union of $F_j$ over all $j \in J$. Since $C_2$ generates $V_I$, and $C_2 \subset F \subset V_I$, this implies that $x_0$ lies in the relative interior of $F$.
\qed 

Lemma \ref{nonor_4.2} allows us to speak of the largest flat of $A$ containing $x_0$ in its interior. There is a particular case where the largest flat of $A$ containing $x_0$ in its interior is easily described :
\begin{lemma}
Take  $x_0$ in  $E$ and  $y \in E^*$ in  $\mathcal{L}(x_0)$. Then any flat of $A$ containing $x_0$ in its interior is contained in $\mathcal{L}(y)$. In particular, if $x_0$ lies in the relative interior of $\mathcal{L}(y)$,  the largest flat of $A$ containing $x_0$ in its interior is   $\mathcal{L}(y)$.
\end{lemma}
\proof
Take a flat $F$ containing $x_0$ in its interior. Then by lemma \ref{interieur_face} 
$\mathcal{F}(F)= \mathcal{L}(x_0)$. Since $y \in  \mathcal{L}(x_0)$, we have $y \in  \mathcal{F}(F)$, that is, for all $x \in F$, $A(x)+B(y) =\langle y,x\rangle  $, i.e $x \in \mathcal{L}(y)$. Hence  $F \subset \mathcal{L}(y)$.
\qed

\section{What we need to know about non-orientable surfaces}
Assume $M$ is non-orientable.
Let $\pi \co M_o \longrightarrow M$ be the orientation cover of
$M$. Then $M_o$ is an orientable surface endowed with a
fixed-point free, orientation-reversing involution $I$. Let
$I_{\ast}$ be the  involution of $H_1 (M_o,\R)$ induced by $I$,
and let $E_1$  (resp. $E_{-1}$) be the eigenspace of $I_{\ast}$ for
the eigenvalue $1$ (resp.$-1$). Likewise, let $I^{\ast}$ be the  involution of $H^1 (M_o,\R)$ induced by $I$,
and let $E_1$  (resp. $E_{-1}$) be the eigenspace of $I^{\ast}$ for
the eigenvalue $1$ (resp.$-1$).
 We have (\cite{nonor},  2.2)
	\[\ker \pi_{\ast}=E_{-1} \subset H_1 (M_o,\R)  \mbox{ and }  \pi^{\ast}\left(H^1(M,\R)\right)=E_{1} \subset H^1 (M_o,\R)
\]
Let 
\begin{itemize}
	\item $T\pi$ denote the tangent map of $\pi$
	\item $L'$ denote the Lagrangian $L \circ T\pi$ on $TM_o$
	\item $\alpha_o$ and $\beta_o$ denote the $\alpha$ and $\beta$-functions, respectively, of $L'$.
\end{itemize}
Likewise we denote with primes the Aubry and Mather sets of $L'$.
Proposition 4 of \cite{Fathi98} says that 
	\[\Azero' = \pi\inv (\Azero), \tildeAzero'=T\pi\inv \left(\tildeAzero\right).
\]
\begin{lemma}\label{Ietoile}
We have
\begin{eqnarray*}
\forall c \in H^1 (M_o,\R),\  \alpha_o (I^* c) & = & \alpha_o ( c)\\	
\forall h \in H_1 (M_o,\R),\  \beta_o (I_* h) & = & \beta_o ( h)
\end{eqnarray*}
\end{lemma}
\proof
Take 
\begin{itemize}
	\item $c \in H^1 (M_o,\R)$
	\item a smooth one-form $\omega$ on $M_o$ such that  $\left[\omega\right]=c$
	\item an $I^*c$-minimizing measure $\mu$.
\end{itemize}
We have
	\[-\alpha_o( I^* c) = \int_{TM_o}\left(L'-I^* \omega \right)d \mu = \int_{TM_o}\left(L'- \omega \right)d I_* \mu \geq -\alpha_o ( c)
\]
where the second equality owes to the $I$-invariance of $L'$. So $ \alpha_o (I^* c) \leq \alpha_o ( c)$, whence $\alpha_o ( c)= \alpha_o (I^* I^* c) \leq \alpha_o ( I^* c)$, which proves the first statement of the lemma.

Now take $h \in H_1 (M_o,\R)$ and an $h$-minimizing measure $\mu$. We have $\left[I_* \mu \right]= I_* h$ thus
	\[\beta_o (I_* h) \leq  \int_{TM_o} L' d I_* \mu =   \int_{TM_o} L' d  \mu = \beta_o ( h)
\]
whence $\beta_o ( h)=\beta_o (I_* I_*  h) \leq \beta_o (I_* h)$, which ends the proof of the lemma.
\qed

\begin{lemma}\label{alpha_revetement}
For all $c \in H^1 (M,\R)$, $\alpha (c)= \alpha_o (\pi^* c)$.
\end{lemma}
\proof

Take $c \in H^1 (M,\R)$ and a smooth one-form $\omega$ on $M$ such that $\left[\omega\right]=c$. Then the lifted Lagrangian corresponding to $L-\omega$ is $L'-\pi^* \omega$. By \cite{Fathi98}, Proposition 4, $L-\omega$ and $L'-\pi^* \omega$ share the same critical value, that is, $\alpha (c)= \alpha_o (\pi^* c)$.

\qed

\begin{lemma}\label{beta_revetement}
For all $h \in E_1 \subset H_1 (M_o,\R)$, $\beta_o (h)= \beta (\pi_* h)$, and if $\mu$ is an $h$-minimizing measure, then $\pi_\ast \mu $ is $\pi_* h$-minimizing.
\end{lemma}
\proof
Take
\begin{itemize}
	\item $h \in E_1 \subset H_1 (M_o,\R)$
	\item an $h$-minimizing measure $\mu$
	\item $c \in H^1 (M_o,\R)$  such that $\alpha_o (c) + \beta_o (h) =\langle c,h\rangle  $.
\end{itemize}
Then by Lemma \ref{Ietoile} $\alpha_o (I^* c) + \beta_o (I_* h) =\langle c,h\rangle  $ and $\langle I^* c,I_* h\rangle  =\langle c,h\rangle  $ since $I^*$ and $I_*$ are dual of one another. Besides, $I_* h =h$ because $h \in E_1$. Setting $c_1 := 2\inv (c+I^* c)$, we have 
	\[ \alpha_o (c_1) \leq \frac{1}{2}\left(\alpha_o (c)+\alpha_o (I^* c)\right)= \alpha_o (c)
\]
by convexity of $\alpha$, but on the other hand
	\[\frac{1}{2}\left(\alpha_o (c)+\alpha_o (I^* c)\right) +\beta_o ( h) = \langle c_1,h\rangle   \leq \alpha_o (c_1)+\beta_o ( h)
\]
whence 
	\[ \alpha_o (c_1) = \frac{1}{2}\left(\alpha_o (c)+\alpha_o (I^* c)\right)= \alpha_o (c)
\]
and 
	\[\langle c_1,h\rangle   = \alpha_o (c_1)+\beta_o ( h).
\]
Since $c_1 \in E_1 \subset H^1 (M_o,\R)$, there exists $c_2$ in $H_1 (M,\R)$ such that $\pi^* c_2= c_1$. By lemma \ref{alpha_revetement}
$\alpha_o (c_1)=\alpha (c_2)$ so 
\begin{eqnarray*}
\alpha (c_2)+	\int_{TM_o} L' d  \mu & =  & \langle \pi^* c_2,h\rangle   \mbox{ that is, }\\
\alpha (c_2)+	\int_{TM} L d \pi_*  \mu & =  & \langle  c_2,\pi_* h\rangle  
\end{eqnarray*}
which proves that $\pi_*  \mu$ is $\pi_* h$-minimizing and $\beta_o (h)= \beta (\pi_* h)$.
\qed

\begin{lemma}\label{mathcalLpi_2}
Let $h$ be an element of $H_1 (M_o,\R)$. We have 
	\[\pi^* \left( \mathcal{L}(\pi_* h)\right)= \mathcal{L}(h)\cap E_1 .
\]
\end{lemma}
\proof
Take $c$ in $\mathcal{L}(\pi_* h)$. We have
	\[ \alpha (c) +\beta (\pi_* h ) = \langle  c,\pi_* h\rangle  
\]
whence by lemmas \ref{alpha_revetement}, \ref{beta_revetement}
	\[ \alpha_o (\pi^* c) +\beta_o ( h ) = \langle  \pi^* c, h\rangle  
\]
that is, $\pi^* c \in \mathcal{L}(h)$. Therefore 
	\[\pi^* \left( \mathcal{L}(\pi_* h)\right)\subset \mathcal{L}(h)\cap E_1. 
\]
Now take $c \in \mathcal{L}(h)\cap E_1$.
Since $c \in E_1$, there exists $c_1 \in H^1 (M,\R)$ such that $\pi^* c_1 =c$. We have 
	\[ \alpha_o ( c) +\beta_o ( h ) = \langle   c, h\rangle   \mbox{ whence } \alpha (c_1) +\beta (\pi_* h ) = \langle  c_1,\pi_* h\rangle  
\]
so $c \in \pi^* \left( \mathcal{L}(\pi_* h)\right)$, which concludes the proof of the lemma.
\qed

\begin{lemma}\label{mathcalLpi}
Let $h$ be an element of $H_1 (M,\R)$, and $h'$ be an element of $E_1 \subset H_1 (M_o,\R)$ be such that $\pi_* h' = h$. We have 
	\[\mathcal{F}(\mathcal{L}(h))= \pi_{\ast}\left(\mathcal{F}(\mathcal{L}(h'))\cap E_1 \right)
\]
\end{lemma}
\proof
Let
\begin{itemize}
	\item $h_1$ be an element of $\pi_{\ast}\left(\mathcal{F}(\mathcal{L}(h'))\cap E_1 \right)$
	\item $h_2$ be an element of $\mathcal{F}(\mathcal{L}(h'))\cap E_1$ such that $\pi_* (h_2)=h_1$
	\item $c$ be an element of $\mathcal{L}(h)$.
\end{itemize}
By Lemma \ref{mathcalLpi_2} $\pi^{\ast}c \in \mathcal{L}(h')$ so 
\[
\alpha_o (\pi^{\ast}c) + \beta_o (h_2)= \left\langle \pi^{\ast}c,h_2 \right\rangle
\]
whence by Lemmas \ref{alpha_revetement},  \ref{beta_revetement}
	\[ \alpha (c) + \beta (h_1)= \left\langle c,h_1 \right\rangle
\]
thus $h_1 \in \mathcal{F}(\mathcal{L}(h))$, hence 
\[\mathcal{F}(\mathcal{L}(h)) \supset \pi_{\ast}\left(\mathcal{F}(\mathcal{L}(h'))\cap E_1 \right).
\]
Conversely, let
\begin{itemize}
	\item $h_1$ be an element of $\mathcal{F}(\mathcal{L}(h))$
	\item $h_2$ be an element of $E_1 \subset H_1(M_o, \R)$ such that $\pi_* (h_2)=h_1$
	\item $c'$ be an element of $\mathcal{L}(h')$.
\end{itemize}	
Setting $c_2 := 2\inv (c'+I^* c')$, we see, as in the proof of Lemma \ref{beta_revetement}, that $c_2 \in \mathcal{L}(h')\cap E_1 $. By Lemma \ref{mathcalLpi_2}, since $c_2 \in \mathcal{L}(h')\cap E_1 $, there exists $c_1$ in 
$\mathcal{L}(h)$ such that $\pi^* c_1= c_2$.

Since $h_1 \in \mathcal{F}(\mathcal{L}(h))$ we have 
	\[ \alpha (c_1) + \beta (h_1)= \left\langle c,h_1 \right\rangle
\]
thus, by Lemma \ref{Ietoile}	
\[\alpha_o (c_2) + \beta_o (h_2)= \left\langle c_2,h_2 \right\rangle .
\]
Hence the two inequalities 
\begin{eqnarray*}
\alpha_o (c') + \beta_o (h_2)& \geq & \left\langle c',h_2 \right\rangle	\\
\alpha_o (I^* c') + \beta_o (h_2)& \geq & \left\langle c',h_2 \right\rangle
\end{eqnarray*}
sum to an equality, so both inequalities are equalities.
Therefore $h_2 \in \mathcal{F}(\mathcal{L}(h'))$, so 
\[\mathcal{F}(\mathcal{L}(h)) \subset \pi_{\ast}\left(\mathcal{F}(\mathcal{L}(h'))\cap E_1 \right).
\]
\qed

\section{Faces and flats of $\beta$}
\label{paragraphe_rationnel}
\subsection{Radial flats of $\beta$}
Recall that $R_h$ is the greatest radial flat of $\beta$ containing the homology class $h$. 
\begin{lemma}\label{Rh_Fh}
For any nonzero $h \in H_1(M,\R)$, for any $t$ such that $th \in R_h$, we have $\mathcal{L}(h)\subset \mathcal{L}(th)$  (in particular, for any non-zero $t$ such that $th \in R_h$, we have $\mathcal{L}(h)= \mathcal{L}(th)$). Consequently,
	\[R_h \subset \mathcal{F}(\mathcal{L}(h)).
\]

\end{lemma}
\proof
Take $t \in \R$ such that $th \in R_h$. By definition of $R_h$ there exists $c \in H^1(M,\R)$ such that 
\begin{eqnarray*}
\alpha(c)+\beta(h)&=&\left\langle c,h\right\rangle\\
\alpha(c)+\beta(th)&=&\left\langle c,th\right\rangle.	
\end{eqnarray*}
The first equality says that $c \in \mathcal{L}(h)$. Take $c'\in \mathcal{L}(h)$. Let us show that $c'\in \mathcal{L}(th)$, which proves that $\mathcal{L}(h) \subset \mathcal{L}(th)$, whence $th \in \mathcal{F}(\mathcal{L}(h))$.

Since $L$ is autonomous,  by \cite{Carneiro} $\alpha(c')=\alpha(c)$ and $\left\langle c',h\right\rangle =\left\langle c,h\right\rangle$.
So 
	\[ \alpha(c')+\beta(th)=\alpha(c)+\beta(th) =\left\langle c,th\right\rangle= \left\langle c',th\right\rangle
\]
that is, $c'\in \mathcal{L}{th}$.
\qed
\begin{lemma}\label{R_h_nonsing}
If  the homology class $h$ is non-singular, then for any $t$ such that   $th \in R_h$,  $th$ is non-singular.
\end{lemma}
Take a non-singular  $h$  and  take  $t$ such that $th \in R_h$. Suppose $th$ is singular, that is, $ \mathcal{M}(\mathcal{L}(th))$ contains a fixed point. The homology of the measure carried by the fixed point is zero, so there exists a cohomology class $c$ such that 
\begin{eqnarray*}
\alpha(c) +\beta(0) & = & \langle c, 0\rangle \\
\alpha(c) +\beta(th) & = & \langle c, th\rangle.
\end{eqnarray*}
Thus there exists a face of $\beta$ containing zero and $th$, since $0 \in \R h$ this means that $0 \in R_{th}$. By the previous lemma this implies 
$$ 0 \in \mathcal{F}(\mathcal{L}(th))$$
but the same lemma says 
$$  \mathcal{F}(\mathcal{L}(th)) \subset  \mathcal{F}(\mathcal{L}(h))$$
whence 
$$ 0 \in \mathcal{F}(\mathcal{L}(h))$$
thus 
$\mathcal{L}(h)\subset \mathcal{L}(0)$
in particular
$$\mathcal{M}(\mathcal{L}(0))  \subset \mathcal{M}(\mathcal{L}(h))$$ 
but $\mathcal{M}(\mathcal{L}(0)) $ contains a fixed point, so $h$ is singular, which is a contradiction.
\qed

 \begin{lemma}\label{radial}
 Let $L$ be a Tonelli Lagrangian on a closed manifold $M$.
 Let $h \in H_1 (M,\R)$ be a  nonsingular homology class. Assume $R(h)=\left[ t_1 h, t_2 h \right]$. Then there exists a sequence of real numbers $t_n$ such that $t_n <t_1$ for all $n$, $t_n$ converges to $t_1$, $t_n h$ is non-singular and 
 $R(t_n h)= \{ t_n h\}$  for all $n$.
 \end{lemma}
 \proof
 The map 
 $$
 \begin{array}{rcl}
\beta_h \co  \R_+ ^* & \longrightarrow & \R \\
 t & \longmapsto & \beta (th)
 \end{array}
 $$
 is convex, superlinear, and $C^1$ (see \cite{Carneiro}). Let $\alpha_h $ be its Fenchel dual, we have 
 $$
 \forall t \in  \R_+ ^* , \   \beta_h (t) +  \alpha_h \left( \beta_h ' (t) \right) = t.  \beta_h ' (t).
 $$
 Let $E$ be the subset of $t \in   \R_+ ^*$  such that $R(th)$ contains properly $\{th\}$. 
 The connected components of $E$ are intervals with non-empty interior, hence $E$ has at most countably many connected components. The derivative of  $\beta_h $ is constant on each connected component of  $E$, hence $ \beta_h ' (E)$ is at most countable.  So the complement in $\R_+ ^*$ of  $ \beta_h ' (E)$  is dense in  $\R_+ ^*$.
 Take a sequence $t_n$ such that   for all $n$, $t_n <t_1$,  
 $ \beta_h ' (t_n) \notin  \beta_h ' (E)$, and   $ \beta_h ' (t_n)$ converges to  $ \beta_h ' (1)$. Then,   since $ \beta_h ' (t_n) \notin  \beta_h ' (E)$, $R(t_n h) =  \{ t_n h\}$  for all $n$. We have 
 $$
\forall n \in  \N , \   \beta_h (t_n ) +  \alpha_h \left( \beta_h ' (t_n) \right) = t_n.  \beta_h ' (t_n)
 $$  
 so by superlinearity of $\beta_h$, the sequence $t_n$ is bounded, and for any of limit-point $t$ of   the sequence $t_n$, we have
 $$
 \beta_h (t ) +  \alpha_h \left( \beta_h ' (1) \right) = t.  \beta_h ' (1)
 $$    
 that is, $t \in \left[ t_1 , t_2  \right]$. Since $t_n <t_1$, we have $t=t_1$.  The set of non-singular homology classes is open in $H_1 (M,\R)$, so for $n$ large enough  $t_n h$ is non-singular. The lemma is proved.
 \qed

Now we look at some consequences of Proposition \ref{rational}.  Assume $h$ is  1-irrational.
Then for all $t$ such that $th \in R_h$, $th$ is also 1-irrational. Furthermore $R_h$ is contained in a face of $\beta$, so Mather's Graph Theorem and Proposition \ref{rational} combine to say that $\tilde{\mathcal{M}}(R_h)$ is a union of pairwise disjoint periodic orbits $\gamma_i$, $i \in I$ where $I$ is some set, not necessarily finite. We denote by $V(R_h)$ the linear subspace of $H_1(M,\R)$ generated by $\left[\gamma_i\right]$, $i \in I$. Since the $\gamma_i$ are pairwise disjoint, there exist homology classes $h_1,\ldots h_k$ with $k \leq 3/2(b_1(M) -2)$, such that $\forall i \in I, \  \exists j=1,\ldots k,  \left[\gamma_i\right]=h_j$.  

Let $T_i$ be the least period of $\gamma_i$. Then the invariant measure $\mu_i$ supported by $\gamma_i$ has homology $T_i\inv \left[\gamma_i\right]$. By Lemma \ref{nonor_4.2} the convex hull $C$ of $T_i\inv \left[\gamma_i\right]$, $i \in I$, is a flat of $\beta$ containing $th$ in its interior whenever $th$ is contained in the  interior of $R_h$. 
\subsection{Faces of $\beta$}
The following lemma is a rewriting of Lemma 12 of \cite{gafa}.
\begin{lemma}\label{face_rationnelle}
Let $F$ be a flat of $\beta$.  Assume $F$ contains a $1$-irrational homology class $h_0$ in its interior. Then for all $h \in F$, for all $h$-minimizing measure $\mu$, the support of $\mu$ consists of periodic orbits, or fixed points.
\end{lemma}
\proof
Let $c$ be a cohomology class  such that for all $h \in F$, $\alpha(c)+\beta(h)=\langle c,h\rangle  $. Take $h$ in $F$. Since $h_0$ lies in the  interior of $F$, there exist $h'$ in $F$ and $0<t<1$ such that $h_0=th+(1-t)h'$. Take an $h$-minimizing (resp. $h'$-minimizing)  measure $\mu$ (resp. $\mu'$), so we have $\beta(h)=\int Ld \mu$,  $\beta(h')=\int Ld \mu'$. Thus  
\begin{eqnarray*}
\alpha(c)+\int Ld \mu &=& \langle c,h\rangle  	\\
\alpha(c)+\int Ld \mu' &=& \langle c,h'\rangle  
\end{eqnarray*}
so
	\[\alpha(c)+\int Ld(t\mu+ (1-t)\mu')= \langle c,th+(1-t)h'\rangle  =\langle c,h_0\rangle  
\]
that is, the probability measure $t\mu+ (1-t)\mu'$ is $h_0$-minimizing. Since $h_0$ is $1$-irrational, Proposition \ref{rational} implies that the support of $t\mu+ (1-t)\mu'$, hence that of $\mu$, consists of periodic orbits, or fixed points.
\qed

Here is a version of Theorem 6.1 of \cite{nonor}  for  general Lagrangians. 
\begin{theorem}\label{thm_gafa}
Let
\begin{itemize}
	\item $M$ be a closed oriented surface
	\item $L$ be a Tonelli Lagrangian on $M$
	\item $h_0$ be  a 1-irrational, nonsingular  homology class of $M$
	\item $V_0$ be $V(R_{h_0})$
	\item $h$ be an element of $V^{\perp}_{0}$.
\end{itemize}
Then there exist $t(h_0,h) \in \R$ and $s(h_0,h) >0$ such that the segment \newline
$\left[h_0,t(h_0,h)h_0 + s(h_0,h)h\right]$ is contained in a face of $\beta$.
\end{theorem}
\proof
We use the notation of Paragraph \ref{paragraphe_rationnel}.

\textbf{First case : }$h \in V_0$. 
Take $t_0$ such that $t_0 h_0$ lies in the relative interior of $R_{h_0}$.  Then  $t_0 h_0$ lies in the relative interior of the convex hull $C$ of
$T_i\inv \left[\gamma_i\right]$, $i \in I$, so there exists a finite subset of $I$, say $\left\{1,\ldots n\right\}$, and $\lambda_1, \ldots \lambda_n$ in $\left] 0,1 \right[$ such that 
\begin{itemize}
	\item $\lambda_1 + \ldots + \lambda_n =1$
	\item $t_0 h_0 = \lambda_1 T_1\inv \left[\gamma_1\right] +\ldots + \lambda_n T_n\inv \left[\gamma_n\right]$
	\item $\left[\gamma_1\right] ,\ldots , \left[\gamma_n\right]$ generate $V_0$.
\end{itemize}
On the other hand, since $h \in V_0$, there exist real numbers $x_1,\ldots x_n$ such that 
$h = x_1 T_1\inv \left[\gamma_1\right] +\ldots + x_n T_n\inv \left[\gamma_n\right]$. Take $\epsilon > 0$ such that $\forall i=1,\ldots n,\  \epsilon x_i + \lambda_i >0$.  Then $(\epsilon x_1 + \lambda_1 +\ldots +\epsilon x_k + \lambda_k)\inv (t_0 h_0 +\epsilon h)$ lies in the relative interior of $C$. Thus there exists a face of $\beta$ containing $h_0$ and $th_0 +sh$, where
	\[ t:= \frac{1}{\sum^{k}_{1}\epsilon x_i + \lambda_i}, \  s := \frac{\epsilon}{t}.
\]

\textbf{Second  case : }$h \notin V_0$, that is, $h \in V^{\perp}_{0} \setminus V_0$. 
\begin{remark}
In that case the dimension of $V_0$ must be less than $b_1(M)/2$.
\end{remark}
For any $n\in \N^{\ast}$ let us denote $h_n := h_0 +n\inv h$. 
Let $\mu_n$ be an $h_n$-minimizing measure.
  For each $i\in I$ let $V_i$ be the neighborhood of $(\gamma_i,\dot\gamma_i)$ given by Lemma \ref{key_o}. Let $V$ be the union over $i\in I$ of the $V_i$. Be sure to take the $V_i$ small enough so $V$ is a disjoint union of annuli.

First let us prove that $V \cap \spt (\mu_n )$ is $\phi_t$-invariant and consists of periodic orbits homotopic to some or all of the $\gamma_i$.
Indeed by Lemma \ref{key_o} a minimizing orbit $\gamma$  that enters $V$ is either
\begin{itemize}
    \item asymptotic to one of the $\gamma_i$
  \item or homotopic to one of the $\gamma_i$
  \item or cuts one of the  $\gamma_i$ with constant sign.  
\end{itemize}
The first case is impossible since $\gamma$ is contained in the support of an invariant measure (see Lemma 5.5 of \cite{nonor}). The third case is impossible since it would imply $\inter(\left[\mu_n\right], \left[\gamma_i\right]) \neq 0$, which contradicts $h \in V^{\perp}_{0}$. So $V \cap \spt (\mu_n )$ is $\phi_t$-invariant and consists of periodic orbits homotopic to some  $\gamma_i$.

Now let us show that for $n$ large enough, $0< \mu_n (V)< 1$. Note that any limit point, in the weak$^\ast$ topology, of the sequence 
$\mu_{n}$ is an $h_0$-minimizing measure, hence supported in $V$, so $\mu_n (V)$ tends to $1$. On the other hand, if $\mu_n (V)$ were $1$, then $\mu_n $ would be supported in $V$.   By the Graph Theorem any minimizing measure supported inside $V$ may be viewed as an invariant measure of a vector field defined in $p(V)$. But $p(V)$ is a union of annuli, hence by the Poincar\'e-Bendixon Theorem any
 minimizing measure supported inside $V$ is supported on fixed points, or periodic orbits homotopic to some  $\gamma_i$. Note that fixed points are ruled out by the nonsingularity of $h$, which implies that $h_n$ is non-singular for $n$ large enough. In particular if $\mu_n (V)=1$, $\left[\mu_n\right] \in V_0$, which is a contradiction. So $0< \mu_n (V)< 1$ and we may set, for any Borelian subset $A$ of $TM$, 
\begin{eqnarray*}
\mu_{n,1} (A) &:=&  \frac{\mu_n(A\cap V)}{\mu_n(V)} \\
\mu_{n,2} (A) &:=&  \frac{\mu_n(A\setminus V)}{\mu_n(TM\setminus V)}\\
\lambda_n &:=& \mu_n(V).
\end{eqnarray*}
Then $\mu_{n,1}$ and $\mu_{n,2}$ are two  probability measures on $TM$. They  are invariant by the Lagrangian  flow because $V\cap \spt (\mu_n )$, as well as its complement in $ \spt (\mu_n )$, is $\phi_t$-invariant. There exists a face of $\beta$ containing $\left[\mu_{n,1}\right]$ and $\left[\mu_{n,2}\right]$ because 
	\[ \mu_1 = \lambda_n \mu_{n,1} + (1-\lambda_n) \mu_{n,2}.
\]
 Let $\mu_{0,1}$ and $\mu_{0,2}$ be limit points, in the weak$^\ast$ topology, of the sequences $\mu_{n,1}$ and $\mu_{n,2}$. Then $\mu_{0,1}$ is an $h_0$-minimizing measure, and there exists a face of $\beta$ containing $h_0 = \left[\mu_{0,1}\right]$ and $\left[\mu_{0,2}\right]$. 

Now we prove that  $\left[\mu_{0,2}\right] \notin V_0$. Assume to the contrary. Then, as in the first case, the face $F_0$ containing $h_0$ and $\left[\mu_{0,2}\right]$ contains $th_0$ in its interior for some $t$ such that $th_0$ lies in $R_0$. Take $\lambda \in \left]0,1\right[$ and $h'$ in $F_0$ such that 
	\[ th_0 = \lambda \left[\mu_{0,2}\right] +(1-\lambda) h'.
\]
Take an $h'$-minimizing measure $\mu'$. Then $\lambda \mu_{0,2} +(1-\lambda) \mu'$ is a $th_0$-minimizing measure, hence is supported inside $V$, which is impossible since $\mu_{n,2}$ is supported outside $V$ for all $n$. 

Thus there exist $v \in V_0$ and $x \neq 0$ such that $\left[\mu_{0,2}\right]=v+xh$. Take $t_0$ such that  $t_0 h_0$ lies in the relative interior of the convex hull $C$ of $T^{-1}_{i}\left[\gamma_i\right], i \in I$, so there exists a positive $\epsilon$ such that $t_0 h_0 -\epsilon v$ lies in the relative interior of $C$.  Lemma \ref{extension_plat} now says that there is a face of $\beta$ containing $h_0$, $t_0 h_0 -\epsilon v$ and $\left[\mu_{0,2}\right]=v+xh$. Such a face must contain $h_0$ and 
	\[\frac{1}{1+\epsilon}\left(t_0 h_0 -\epsilon v \right)+ \left(1- \frac{1}{1+\epsilon}\right)\left(v+xh \right)= 
	\frac{t_0}{1+\epsilon}h_0 +\frac{x \epsilon }{1+\epsilon}h,
\]
Now set 
	\[ t(h_0,h) := \frac{t_0}{1+\epsilon},\  s(h_0,h) := \frac{x \epsilon }{1+\epsilon}
\]
and the theorem is proved.
\qed

{\small

\bigskip

\noindent

D\'epartement de Math\'ematiques, Universit\'e Montpellier 2, France\\
e-mail : massart@math.univ-montp2.fr
}

\end{document}